\documentclass[1p]{elsarticle}

\usepackage{amsmath,amsthm,amssymb,amsbsy,epsfig,graphicx,color,multicol,subfigure}
\usepackage{array,calc}
\usepackage[enableskew]{youngtab}
\usepackage{rotating}
\usepackage{young,epic}
\usepackage{a4wide,bm}
\usepackage{url}
\usepackage{hyperref}

\def\nn{\nonumber}

\def\ba{\begin{eqnarray}}
\def\ea{\end{eqnarray}}

\newcommand{\nc}{\newcommand}
\nc{\be}{\begin{equation}} \nc{\ee}{\end{equation}}
\nc{\bea}{\begin{eqnarray}} \nc{\eea}{\end{eqnarray}}
\nc{\ben}{\begin{eqnarray*}}
\nc{\een}{\end{eqnarray*}}
\nc{\mmatrix}[1]{\begin{matrix} #1 \end{matrix}}
\nc{\dt}[1]{\left| \begin{matrix} #1 \end{matrix}\right|}
\nc{\x}{\mathbf x}
\nc{\y}{\mathbf y}
\nc{\z}{\mathbf z}
\nc{\av}{\mathbf a}
\newenvironment{invsquarecells}[2]
{\setlength\celldim{#2}%
\settoheight\fontheight{A}%
\setlength\extraheight{\celldim - \fontheight}%
\setcounter{sqcolumns}{#1 - 1}%
\begin{tabular}{|S|*{\value{sqcolumns}}{Z|}}}
{\end{tabular}}
\nc{\disp}{\displaystyle} \nc{\ade}{\mbox{$A$-$D$-$E$}}
\nc{\calN}{{\cal N}} \nc{\calC}{{\cal C}} \nc{\calM}{{\cal M}}
\nc{\calS}{{\cal S}} \nc{\phit}{\hat{\varphi}}
\nc{\chit}{\hat{\chi}} \nc{\hcalN}{\hat{\calN}}
\nc{\hcalS}{\hat{\calS}} \nc{\hS}{\hat{S}}
\nc{\sigmad}{\sigma^\dagger} \nc{\psid}{\psi^\dagger}

\def\d{{\rm d}}
\def\i{{\rm i}}
\def\H{\mathcal{H}}

\def\KL{{\rm KL}}
\def\M{{\rm M}}

\definecolor{IndianRed}{rgb}{0.8,0.36,0.36}
\newtheorem{example}{Example}
\newtheorem{remark}{Remark}
\newtheorem{lemma}{Lemma}
\newtheorem{theorem}{Theorem}

\newtheorem{prop}{Proposition}
\newtheorem{cor}{Corollary}

\newcommand{\crot}[1]{\begin{turn}{45} #1 \end{turn}}
\newlength\celldim \newlength\fontheight \newlength\extraheight
\newcounter{sqcolumns}
\newcolumntype{S}{ @{}
>{\centering \rule[-0.5\extraheight]{0pt}{\fontheight + \extraheight}}
p{\celldim} @{} }
\newcolumntype{Z}{ @{} >{\centering} p{\celldim} @{} }

\newenvironment{squarecells}[2]
{\setlength\celldim{#2}%
\settoheight\fontheight{A}%
\setlength\extraheight{\celldim - \fontheight}%
\setcounter{sqcolumns}{#1 - 1}%
\begin{tabular}{|S|*{\value{sqcolumns}}{Z|}}\hline}
{\end{tabular}}
\newcommand\nl{\tabularnewline\hline}
\DeclareMathOperator{\CT}{CT}
\begin{document}
\begin{frontmatter}

\title{Deformed Kazhdan-Lusztig elements and Macdonald polynomials}
\author[jdg]{Jan de Gier}
\ead{jdgier@unimelb.edu.au}
\address[jdg]{Department of Mathematics and Statistics, The University of Melbourne, VIC 3010, Australia}

\author[al]{Alain Lascoux}
\ead{al@univ-mlv.fr}
\address[al]{Universit\'e Paris-Est, CNRS, Institut Gaspard Monge, 77454 Marne La Vallee, Cedex 2, France}

\author[jdg]{Mark Sorrell}
\ead{msorrell@unimelb.edu.au}

\begin{abstract}
We introduce deformations of Kazhdan-Lusztig elements and  specialised  nonsymmetric Macdonald polynomials, both of which form a distinguished basis of the polynomial representation of a maximal parabolic subalgebra of the Hecke algebra. We give explicit integral formula for these polynomials, and explicitly describe the transition matrices between classes of polynomials. We further develop a combinatorial interpretation of homogeneous evaluations using an expansion in terms of Schubert polynomials in the deformation parameters.
\end{abstract}

\end{frontmatter}

\section{Introduction}
We discuss two important bases in the Iwahori-Hecke algebra $\H$: the Kazhdan-Lusztig basis \cite{KL79} and Young's seminormal, or orthogonal basis \cite{Young,Hoefsmit}. While the former admits a deep geometric interpretation, the latter is purely algebraic and combinatorial. We study these two bases for the maximal parabolic case in a very explicit way using the polynomial representation of $\H$. In this representation, the parabolic KL basis \cite{Dheod87} gives rise to KL elements, while the orthogonal basis elements are given by specialised non-symmetric Macdonald polynomials \cite{Cher92,BernardGHP,Macd96}. We review their construction and for both cases give elegant explicit formulas using factorisations in terms of Baxterised operators. We further introduce natural deformations which interpolate bewteen classes of KL elements and specialised Macdonald polynomials, as well as provide explicit multi-dimensional integral expressions for the polynomial basis elements.

The motivation for studying the transition matrices between these two bases arises from physics as well as from combinatorics. Both distinguished bases have important applications in physics when used in finite dimensional polynomial representation modules of the Hecke algebra. For example, they are related to quantum incompressible states \cite{Pasq05}, which in the limit to Jack polynomials contain model groundstate wave functions of (fractional) quantum Hall systems \cite{BerneH08}. They are also used to describe polynomial solutions to the $q$-deformed Knizhnik-Zamolodchikov (qKZ) equation \cite{DFZJ05,KasatPasq06,KasatTake07}. The latter are relevant for, for example, critical bond percolation in an inhomogeneous system \cite{DFZJ04,DFZJ07,DFZJ07b,ZJ07,GierPS,Cant09,GierNP10}. The investigation of the transition matrix between the KL and the orthogonal basis was inspired by this latter application, where naturally arising sums over KL basis elements were observed to equal just one single non-symmetric Macdonald polynomial. Here we prove this fact. In the context of physics applications, explicit expressions are also needed for the homogeneous evaluations of these polynomials, and we study these with the theory of Schubert polynomials as a useful tool. 

From a different point of view, in a recent preprint \cite{Blasiak11} similar transition matrices between inequivalent KL bases of Specht modules are studied, as well as their relation to canonical bases of corresponding quantum group modules using quantum Schur-Weyl duality. 

There is a rich combinatorial content in this theory, and we proceed to derive combinatorial rules for certain transition matrices between basis elements. In addition, the homogeneous evaluations we consider give rise to enumeration formulas for combinatorial objects such as fully packed loop configurations, alternating sign matrices and punctured symmetric plane partitions, see e.g. \cite{Bressoud}. Using the theory of Schubert polynomials as a computational tool we obtain, as a corollary, natural combinatorial interpretations of the deformation parameters.

\section{Polynomial representation of the Hecke algebra}
The Hecke algebra $\H$ of the symmetric group $W=S_n$ generated by the simple reflections $s_i$, is the algebra over $\mathbb{Q}[t,t^{-1}]$ defined in terms of generators $T_i\equiv T_{s_i}$, $i=1,\ldots,n-1$, and relations
\begin{align}
&(T_i-t)(T_i+t^{-1}) = 0,\qquad T_iT_j = T_jT_i \quad\forall\, i,j:\, |i-j| >1,\nonumber\\
&T_iT_{i+1}T_i = T_{i+1}T_iT_{i+1}. \label{Hdef_typeA}
\end{align}
The Baxterised element $T_i(u) \in \mathbb{C}\H$,  parametrised by the complex number $u$, is defined by
\[
T_i(u) = T_i + \frac{t^{-u}}{[u]_t},
\]
where the notation $[u]_t$ stands for the usual t-number
\[
[u]_t = \frac{t^u-t^{-u}}{t-t^{-1}}.
\]
The notation $[u]$ will always refer to base $t$. The element $T_i(u)$ is constructed to satisfy the Yang-Baxter equation,
\be
T_i(u)T_{i+1}(u+v)T_i(v) = T_{i+1}(v)T_i(u+v)T_{i+1}(u),
\label{eq:ybe}
\ee
which will be important later on.

The standard projectors are obtained by specialising $u$:
\[
T_i(1) = T_i+t^{-1},\qquad T_i(-1) = T_i-t^{}.
\]

The Hecke algebra has a multivariate polynomial representation, which is conveniently expressed in terms of the divided difference operator $\partial_i$, defined by
\be
\label{pdiff}
\partial_i f = \frac{f-s_i f}{z_i-z_{i+1}},
\ee
where $s_i$ denotes the transposition $z_i\leftrightarrow z_{i+1}$. The projector $T_i(1)$ then induces the operator $\nabla_i$.  With abuse of notation: 
\be
T_i(1) = T_i +t^{-1} = \nabla_i = (t z_i -t^{-1}z_{i+1})\partial_i := \frac{t z_i-t^{-1}z_{i+1}}{z_i-z_{i+1}} (1-s_i), 
\ee
which commutes with multiplication with functions symmetric in $z_i,z_{i+1}$, and acts on $1$ and $z_i$ as
\be
\nabla_i : \left\{ 
\renewcommand{\arraystretch}{1.4}
\begin{array}{rcl}
z_i & \mapsto & t z_{i} - t^{-1} z_{i+1}\\ 
1 & \mapsto & 0
\end{array}\right.
\ee
Representations of the Hecke algebra are labeled by partitions. In this paper we will restrict ourselves to  representations  corresponding to rectangular partitions with two rows or two columns  of length $n$ . In such representations, basis elements can be labeled by partitions contained in the maximal staircase partition  $(n-1,n-2,\ldots,1)$ . It is our aim to describe several distinguished polynomial bases.

We will sometimes indicate partitions by their corresponding Yamanouchi word. Recall that a Yamanouchi word $w$ is a word with integer entries such that for every factorisation $w=w' w''$, the right factor $w''$ contains as many or more occurrences of the symbol $i$ than of $i+1$, for all $i \ge 0$. Given any Yamanouchi word, a dual Yamanouchi word may may be obtained by the following procedure: Begin by numbering, from right to left, each integer of the Yamanouchi word that is equal to $0$.  Then repeat for each integer equal to $1$, and so on up to the maximal integer appearing in the Yamanouchi word.  The resultant label from this numbering is also a Yamanouchi word, which we call its dual.
\begin{example}
We begin with the Yamanouchi word $11011000$, and follow the procedure described above. For any $i$, label each occurence of $i$ from right to left, by $0,1,2,\ldots$:
$$ \begin{array}{cccc cccc} 
 1 &1 &0 &1 &1 &0 &0 &0 \\  \noalign{\hrule}
   &  &3 &  &  &2 &1 &0 \\
 3 &2 &  &1 &0 &  &  &  \\ \noalign{\hrule}
 3 &2 &3 &1 &0 &2 &1 &0     
\end{array}
$$
This yields the Yamanouchi word $32310210$ dual to $11011000$.
\end{example}

There is a further simple bijection between Yamanouchi words with two distinct integers, and sub-partitions of a staircase partition.  This is described as follows. Labeling each vertical step with a $1$, and each horizontal step with a $0$, the staircase partition is labeled $1010\dots 10$, which is the Yamanouchi word dual to the word $(n-1)(n-1)\dots 1100$.  Any sub-partition of this staircase will have a lower edge which is a Dyck path, and its labeling in terms of $1$s and $0$s will form a Yamanouchi word, e.g. the Dyck path arising from $(1,1) \subset (3,2,1)$ is 

\[
\raisebox{-30pt}{\includegraphics[height=60pt]{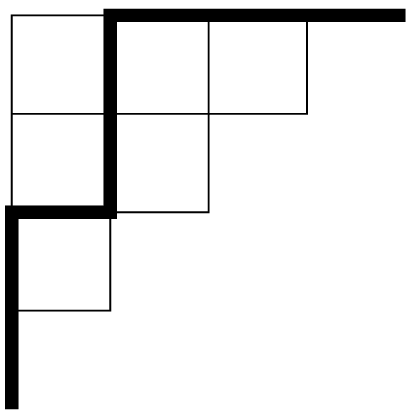}}\quad \rightarrow\quad 11011000.
\]

\subsection{Kazhdan-Lusztig basis from $t$-Vandermonde determinant}

Kazhdan and Lusztig (KL) \cite{KL79} defined a linear basis of the Hecke algebra $\H$, in relation with a fundamental involution of $\H$. 
%%%%%%%%%%%%%%%%%%%%%%%%%%
 Given any polynomial $f(z_1 \dots z_{2n})$, one can use the KL basis to study the module $\H\, f$. In particular, defining the $t$-Vandermonde by
\be
\Delta_t(z_1,\ldots,z_{N})= \prod_{1\leq i<j\leq N} (tz_i-t^{-1}z_{j}).
\label{tVand}
\ee
  and $\Delta\Delta$ to be 
$$ \Delta\Delta = \Delta_t(z_1,\ldots,z_{n}) \Delta_t(z_{n+1},\ldots,z_{2n}),
%\prod_{1\leq i<j\leq n} (tz_i-t^{-1}z_{j})  \prod_{n+1\leq i<j\leq 2n} (tz_i-t^{-1}z_{j})  
$$
then one shows  \cite{KiriL00} that the module $\H\, \Delta\Delta$ is an irreducible representation of $\H$ corresponding to the partition $(n,n)$, and that a subset of the Kazhdan-Lusztig basis furnishes a basis, that we shall still call KL basis, of this space.   

In this special case, the elements of the KL basis can be labeled by partitions contained in the staircase $(n-1,\ldots,1)$, or alternatively, by Yamanouchi words. They furthermore satisfy two special properties. The first is that they satisfy the following vanishing property. Let $w,w'\in \mathbb{N}^{2n}$ be Yamanouchi words, and let $t^{-2w}=\{t^{-2w_1},\ldots,t^{-2w_{2n}}\}$, then
\be
\KL_{w'}(z_1,\ldots,z_{2n})\Big|_{z=t^{-2w}} = 0 \quad \textrm{unless}\quad w=w'.
\ee

The second property is that they can be obtained from $\Delta\Delta$  by applying specific products of $T_i(u)$  \cite{KiriL00,GierP07}. For example, for $n=3$ we have:

\begin{align*}
\KL_{111000} &=    \Delta\Delta= 
\Delta_t(z_1,z_2,z_{3})\Delta_t(z_{4},z_5,z_{6}),\\
\KL_{110100} &= T_3(1)\ \KL_{111000},\\
\KL_{101100} &= T_2(1)T_3(2)\ \KL_{111000},\\
\KL_{110010} &= T_4(1)T_3(2)\ \KL_{111000},\\
\KL_{101010} &= T_2(1)T_4(1)T_3(2)\ \KL_{111000}.
\end{align*}
The order in which the Baxterised operators $T_i(u)$ are applied is determined by a Young diagram. For example, if we graphically denote $T_i(u)$ by a labeled tilted square
\[
T_i(u) = \raisebox{6pt}{
\begin{turn}{-45} 
\begin{squarecells}{1}{1.5em}
\crot{$u$} \nl
\end{squarecells}
\end{turn}},
\]
acting from top to bottom in the $i$th column, then the polynomial

\[
\KL_{1101101000} = T_6(1)T_3(1)T_4(2)T_5(3)\ \KL_{1111100000},
\]
may be graphically depicted  as
\[
\begin{turn}{-45}
\begin{squarecells}{4}{1.5em}
\crot{3} & \crot{1} & \multicolumn{2}{l}{\hphantom{\crot{2}\crot{2}} }\vline \cr\cline{1-2}\cline{4-4}
\crot{2} & \multicolumn{2}{c}{\hphantom{\crot{2}\crot{2}} }\vline \cr\cline{1-1}\cline{3-3}
\crot{1} &  \cr\cline{1-2}
\cr\cline{1-1}
\end{squarecells}
\end{turn}\ .
\]

Rotating these pictures by $\pi/4$, it is clear that we may index the KL basis by labeled partitions that fit inside the staircase, the index $i$ of each operator $T_i(u)$ being determined by the position of the box.

\begin{align*}
\KL_{111000} &=  \KL_{\emptyset} = \begin{squarecells}{2}{1em} 
 & \nl
\cr\cline{1-1}
\end{squarecells}\ ,\quad \big(= \Delta_t(z_1,z_2,z_{3})\Delta_t(z_{4},z_5,z_{6})\big),\\
\KL_{110100} &=  \KL_{(1)}  = T_3(1)\ \KL_{\emptyset} =  
\begin{squarecells}{2}{1em} 
1 & \nl
\cr\cline{1-1}
\end{squarecells}\ ,\\
\KL_{101100} &=   \KL_{(1,1)}    =  T_2(1)\ T_3(2)\ \KL_{\emptyset} =
\begin{squarecells}{2}{1em}
2 & \nl
1 \cr\cline{1-1}
\end{squarecells}\ ,\\
\KL_{110010} &=   \KL_{(2)}    = 
\begin{squarecells}{2}{1em}
2 & 1 \nl
\cr\cline{1-1}
\end{squarecells}\ \\
\KL_{101010} &=   \KL_{(2,1)}    = 
\begin{squarecells}{2}{1em}
2 & 1 \nl
1\cr\cline{1-1}
\end{squarecells}\ .
\end{align*}

The general rule for the labels associated to KL elements given by \cite{KiriL00,GierP07} can be summarised as follows. Let $v_{ij}$ be the label for box $(ij)$ in the inner shape, then the labels are determined by the recursion
\[
v_{ij} = \max\{v_{i+1,j},v_{i,j+1}\} +1,
\]
where $v_{ij}=0$ if $(ij)$ lies outside the inner shape. Note that the value of the labels depend on the inner shape. For example,
\[
\setlength{\unitlength}{1pt}
\raisebox{12pt}{
\begin{picture}(150,150)
\linethickness{2pt}
\put(0,55){\begin{squarecells}{7}{2em}
& & & & & & \nl
& $v_{i,j}$ & & & & \cr\cline{1-6}
& & & & \cr\cline{1-5}
& & & \cr\cline{1-4}
& & \cr\cline{1-3}
& \cr\cline{1-2}
\cr\cline{1-1}
\end{squarecells}}
\put(-1,47){\line(1,0){20}}
\put(20,46){\line(0,1){41}}
\put(19,88){\line(1,0){42}}
\put(61,87){\line(0,1){20}}
\put(60,107){\line(1,0){22}}
\put(81,108){\line(0,1){21}}
\put(82,129){\line(-1,0){83}}
\put(0,129){\line(0,-1){82}}
\end{picture}},
\qquad\qquad
\raisebox{75pt}{$\Longrightarrow$} \qquad\qquad
\raisebox{12pt}{
\begin{picture}(150,150)
\linethickness{2pt}
\put(0,55){\begin{squarecells}{7}{2em}
4 & 3 & 2 & 1 & & & \nl
3 & 2 & 1 & & & \cr\cline{1-6}
2 & & & & \cr\cline{1-5}
1 & & & \cr\cline{1-4}
& & \cr\cline{1-3}
& \cr\cline{1-2}
\cr\cline{1-1}
\end{squarecells}}
\put(-1,47){\line(1,0){20}}
\put(20,46){\line(0,1){41}}
\put(19,88){\line(1,0){42}}
\put(61,87){\line(0,1){20}}
\put(60,107){\line(1,0){22}}
\put(81,108){\line(0,1){21}}
\put(82,129){\line(-1,0){83}}
\put(0,129){\line(0,-1){82}}
\end{picture}}
\]

\subsection{Young basis from $t$-Vandermonde determinant}

There are other bases of the module $\H \Delta\Delta$  that can be obtained as images of  $\Delta\Delta$ under  products of $T_i(u)$, which 
are specified by the choice of an arbitrary vector of parameters called spectral vector. These bases generalise Young's orthonormal basis  of irreducible representations, which corresponds to choosing the vector of contents as spectral vector \cite[section 2]{Hoefsmit,Las01}. 

In our case, one of the bases can be obtained as a specialisation of interpolation Macdonald polynomials  \cite{Knop97,Sahi98}. Indeed, the interpolation polynomial of index $n-1, \ldots 0, n-1, \ldots , 0$   specializes for $q = t^6$ to the product  $\Delta\Delta$ \cite{KasatPasq06}, and the set of homogeneous polynomials $M_v$,  for all  Yamanouchi words which are a permutation of $n-1, \ldots 0, n-1, \ldots , 0$, is a linear basis of $\H\Delta\Delta$. 

Like the KL basis, this basis of specialised Macdonald polynomials may be indexed by labeled partitions $\lambda\subseteq\rho=(n-1,\ldots,1)$.  We shall denote these polynomials $M_\lambda$. In this case the labels do not depend on the partition $\lambda$ as was the case for the KL basis. Rather the labels of the subdiagram are just those inherited from $\rho$, the maximal diagram. For example, if $\lambda=(4,3,1,1)$, then the labels of both the maximal diagram and $\lambda$ are as in the following figure,
\[
\setlength{\unitlength}{1pt}
\raisebox{12pt}{
\begin{picture}(150,150)
\linethickness{2pt}
\put(0,55){\begin{squarecells}{7}{2em}
$\scriptstyle n$ & $\scriptstyle n-1$ & $\scriptstyle n-2$ & $\scriptstyle n-3$ & \ldots & 3 & 2\nl
$\scriptstyle n-1$ & $\scriptstyle n-2$ & $\scriptstyle n-3$ & & & 2 \cr\cline{1-6}
$\scriptstyle n-2$ & & & & \cr\cline{1-5}
$\scriptstyle n-3$ & & & \cr\cline{1-4}
\vdots & & \cr\cline{1-3}
3 & 2 \cr\cline{1-2}
2 \cr\cline{1-1}
\end{squarecells}}
\put(0,48){\line(1,0){21}}
\put(20,48){\line(0,1){42}}
\put(20,89){\line(1,0){42}}
\put(61,89){\line(0,1){21}}
\put(61,109){\line(1,0){21}}
\put(81,108){\line(0,1){21}}
\put(82,129){\line(-1,0){83}}
\put(0,129){\line(0,-1){82}}
\end{picture}}
\]

\begin{example}
For $n=3$ the irreducible module corresponds to
\begin{align*}
\M_{210210} &= \Delta_t(z_1,z_2,z_{3})\Delta_t(z_{4},z_5,z_{6})\\
\M_{212010} &= T_3(3)\ \M_{210210}\\
\M_{221010} &=T_2(2)\ T_3(3)\ \M_{210210}\\
\M_{212100} &=T_4(2)\ T_3(3)\ \M_{210210}\\
\M_{221100} &=T_2(2)\ T_4(2)\ T_3(3)\ \M_{210210}\\
\end{align*}
which we abbreviate by
$$
\begin{array}{lll}
\M_{210210} = \emptyset,\qquad &
\M_{212010} = \begin{squarecells}{1}{1em}
$3$ \nl
\end{squarecells}\ ,\\[2mm]
\M_{221010} = \begin{squarecells}{1}{1em}
$3$ \nl
$2$\cr\cline{1-1}
\end{squarecells}\ , &
\M_{212100} = \begin{squarecells}{2}{1em}
$3$ & $2$ \nl
\end{squarecells}\ , &
\M_{221100} = \begin{squarecells}{2}{1em}
$3$ & $2$ \nl
$2$\cr\cline{1-1}
\end{squarecells}\ .
\end{array}
$$
\end{example}

Given the KL and M bases, it is natural to look at the transition matrices between them.  For example, let $[\rho]=1010\ldots 10$ denote the maximal Yamanouchi word corresponding to the staircase partition $\rho=(n-1,\ldots,1)$. Then $[\tilde{\rho}]=n-1,n-1,\ldots,221100$ is the Yamanouchi word dual to $[\rho]$. The expansion of the maximal basis element $\M_{[\tilde{\rho}]}$ in terms of KL polynomials was conjectured to be given by (equation (5.3) in \cite{GierP07})
\be
M_{[\tilde{\rho}]} =: M_{\rho} = \sum_{\lambda\subseteq \rho}\ \tau^{-n_\lambda}\ \KL_{\lambda},\qquad \tau=-[2].
\label{M2KL}
\ee

Here $n_\lambda$ is defined as the signed sum of boxes between the Dyck path corresponding to the maximal staircase $\rho$, and the Dyck path corresponding to the partition $\lambda$. This is most easily seen via an example:
\begin{example}
Consider $n=6$, with $\lambda=(3,1)$ as shown in Figure \ref{cdef}. In this example we thus have $n_\lambda = 5-4+2 = 3$.

\begin{figure}[here!]
\begin{center}
\resizebox{6cm}{!}{\includegraphics{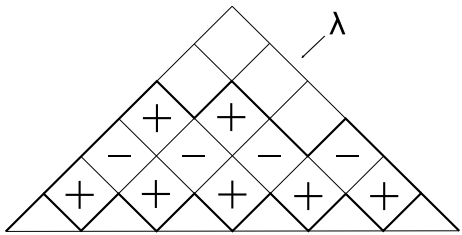}}
\caption{Definition of $n_\lambda$ as the sum of signed boxes. The first row of boxes is labelled with plus signs, the next with minus signs, and so on.
\label{cdef}}
\end{center}
\end{figure}

\end{example}
We will prove a more general result than \eqref{M2KL} in Theorem~\ref{th:expansion}, see Section~\ref{se:expansionM2KL}.

\section{Deformations}

  Let us take a pair of vectors, $w \in\mathbb{N}^{2n}$ and $\langle w\rangle \in \mathbb{C}^{2n} $, called  \textit{spectral vector}, satisfying the property  that, for each  $i$ such that $w_i < w_{i+1}$,  then $\langle w\rangle_{i+1} \neq \langle w\rangle_{i}$ and $\langle s_i w\rangle =s_i \langle w\rangle$, where $s_i$ is a simple transposition. From a starting polynomial $M_{w_0}$ in $z_1, \ldots, z_{2n}$,   one generates other polynomials by the recursion    
\[
M_w \rightarrow M_{s_i w} := T_i\big(\langle w\rangle_{i+1}-\langle w\rangle_{i}\big)M_v.
\]
Thanks to the Yang-Baxter equation \eqref{eq:ybe}, the polynomials $M_{w'}$ are independent of the choice of the reduced decomposition chosen to pass from $w$ to $w'$.

Let us now be more specific and take 

\[
w_0 = (n-1,\ldots, 0, n-1,\ldots, 0),\qquad \langle w_0\rangle = (n-1,n-2-u_1,n-3-u_2,\cdots,-u_{n-1},n-1,n-2,\ldots,2,1,0).
\]

As we did in the case of KL polynomials, let us rather use partitions for indexing. 
  
The starting point is $M_{\emptyset}=M_{w_0}=
\Delta\Delta=
\Delta_t(z_1,\ldots,z_{n})\Delta_t(z_{n+1},\ldots,z_{2n})$.
The rule to obtain  the other Macdonald polynomials, which is equivalent
to the description in terms of spectral vectors given in \cite{Las01}, is as follows. Given an arbitrary partition $\lambda=(\lambda_{n-1},\ldots,\lambda_1)$ with $\lambda_{n-1} \geq \lambda_{n-2} \geq \ldots \geq \lambda_1$, the label of box $(ij)$ of the corresponding diagram is equal to $u_{n-i} + n -i - j+1$, for example,
\[         
\setlength{\unitlength}{1pt}
M_{(\lambda_{n-1},\ldots,\lambda_1)}(u_1,\ldots,u_{n-1}) = \;
\raisebox{-120pt}{ 
\begin{picture}(220,220) 
\setlength{\unitlength}{1pt}
\linethickness{2pt}
\put(0,105){\begin{squarecells}{6}{35pt}
\parbox{30pt}{$\scriptstyle u_{n-1} +$ \\ \mbox{} \hspace{10pt} $\scriptstyle n-1$} & \parbox{30pt}{$\scriptstyle u_{n-1} +$ \\ \mbox{} \hspace{10pt} $\scriptstyle n-2$} & \ldots & \ldots & \ldots & $\scriptstyle u_{n-1}+1$\nl
\vdots & & & &  \cr\cline{1-5}
\vdots & & & \cr\cline{1-4}
\vdots& & \cr\cline{1-3}
$\scriptstyle u_2+2$ & $\scriptstyle u_2+1$ \cr\cline{1-2}
$\scriptstyle u_1+1$ \cr\cline{1-1}
\end{squarecells}}
\put(0,72){\line(1,0){36}}
\put(35,72){\line(0,1){71}}
\put(35,142){\line(1,0){72}}
\put(106,142){\line(0,1){37}}
\put(106,178){\line(1,0){36}}
\put(142,177){\line(0,1){36.9}}
\put(142,213){\line(-1,0){143}}
\put(0,213){\line(0,-1){142}}
\end{picture}}
\]
where we have suppressed the dependence on $z_1,\ldots,z_{2n}$. The undeformed Macdonald polynomials are obtained by setting $u_i=1$.

\begin{example}
The deformed maximal Macdonald polynomial (for $n=3$, $q=t^6$) is given by
\[
\M_{(2,1)}(u_1,u_2)=\M_{221100}(u_1,u_2) = \begin{squarecells}{2}{2.2em}
$\scriptstyle u_2+2$ & $\scriptstyle u_2+1$ \nl
$\scriptstyle u_1+1$\cr\cline{1-1}
\end{squarecells}\  = T_2(u_1+1)T_4(u_2+1)T_3(u_2+2)\ M_{210210}.
\]
\end{example}

\subsection{Integral formulas}
\label{integrals}
Introducing an auxiliary variable $w$, we can separate the variables $z_i$ and $z_{i+1}$ and note the following
\begin{multline}
T_i(1) \frac{\Delta_t(z_i,z_{i+1})}{(w-z_i)(tw-t^{-1}z_{i+1})} =\\\
 \Delta_t(z_i,z_{i+1})\left(\frac{1}{(w-z_i)(w-z_{i+1})} + \frac{1}{(tw - t^{-1}z_i)(tw-t^{-1}z_{i+1})}\right).
\label{Tonphi}
\end{multline}
This observation naturally leads to integral formulas using Bethe Ansatz functions, see e.g. \cite{JimboMiwa} for similar formulas in the context of solutions to the $q$-deformed Knizhnik-Zamolodchikov equations. First we define
\be
\phi_i(w) = \prod_{m=1}^i \frac{1}{w-z_m} \prod_{m=i+1}^{2n}\frac{1}{tw-t^{-1}z_m},
\ee
and introduce the shift operator $S_i$ which acts as
\be
S_i \phi_i = \phi_{i-1},\qquad S_i \phi_j = \phi_j\quad (j\neq i).
\ee
Using \eqref{Tonphi}, we have the following simple action,
\begin{align*}
T_i(v_k+1) \Delta_t(z_1,\ldots,z_N)\phi_{i}(w) &=  \Delta_t(z_1,\ldots,z_N) \Big(S_i + S_{i+1}^{-1} + y_k\Big) \phi_i(w),\\
T_i(v_k+1) \Delta_t(z_1,\ldots,z_N)\phi_{j}(w) &= \frac{[2+v_k]}{[1+v_k]} \Delta_t(z_1,\ldots,z_N)\phi_{j}(w)\qquad (j\neq i),
\end{align*}
where
\[
y_k = -\frac{[v_k]}{[v_k+1]}.
\]
The action of the Hecke generators on the ``wave functions" $\phi_i$ leads to the following result, in which the deformed M polynomials are explicitly given as multiple integrals.

\begin{theorem}
\label{DefMacdIntegral}

Set $\lambda_0=u_0=0$, and let $a_k=\lambda_k+k+1$ for $0\leq k\leq n-1$. Let furthermore $y_k = -\frac{[v_k]}{[v_k+1]}$ where $v_k= u_k+k-\lambda_k$. We have
\begin{align*}
&\M_{(\lambda_{n-1},\ldots,\lambda_1)}(u_1,\ldots,u_{n-1}) \\
&\hphantom{\M_{(\lambda_{n-1},\ldots,\lambda_1)}} = \Delta_t(z_1,\ldots,z_{2n}) \oint  \cdots \oint \frac{1}{(2\pi\i)^n} \Delta(w_n,\ldots,w_1) \Delta_t(w_1,\ldots,w_n) \times\\
&\hphantom{\M_{(\lambda_{n-1},\ldots,\lambda_1)} \Delta_t(z_1,\ldots,z_{2n})} \prod_{m=0}^{n-1} \left(1+y_{m} S_{a_{m}}\right)\ \phi_{a_{m}}(w_{m+1})\hspace{1mm}\d w_1 \ldots \d w_n\\
&\hphantom{\M_{(\lambda_{n-1},\ldots,\lambda_1)}} = \Delta_t(z_1,\ldots,z_{2n}) \oint \cdots \oint \frac{1}{(2\pi\i)^n}\ \Delta(w_n,\ldots,w_1) \Delta_t(w_1,\ldots,w_n) \times\\
& \hphantom{\M_{(\lambda_{n-1},\ldots,\lambda_1)} \Delta_t(z_1,\ldots,z_{2n})} \prod_{m=0}^{n-1} \frac{1}{[v_{m}+1]}\left(\frac{t^{v_{m}+1} w_{m+1} -t^{-v_{m}-1}z_{a_{m}}}{tw_{m+1}-t^{-1}z_{a_{m}}} \right) \phi_{a_{m}}(w_{m+1}) \hspace{1mm}\d w_1 \ldots \d w_n.
\end{align*}
\end{theorem}

\begin{remark}
Recall that the undeformed M polynomials are obtained by setting $u_k=1$ for all values of $k$. If $\lambda$ is a strict partition, then the deformed M polynomials are equal to the KL polynomials when  $v_k=0$, or in other words $u_k=\lambda_k-k$ .
\end{remark}

\begin{proof}
The proof is inductive and given in \ref{se:DefMacdIntegralProof}.
\end{proof}
\bigskip

\subsection{Expansion of M into KL polynomials}
\label{se:expansionM2KL}

Having the integrals of Theorem~\ref{DefMacdIntegral} at our disposal, we now investigate the expansion of M-basis elements into the KL basis.

\begin{theorem}
\label{th:expansion}
Let $\mu=(\mu_{n-1},\ldots,\mu_1) =(n-1,\ldots,1)$ be the staircase partition. The deformed maximal Macdonald polynomial $M_{\mu}(u_1,\ldots,u_{n-1})$ is equal to the sum
\[
M_{\mu}(u_1,\ldots,u_{n-1}) = \sum_{\lambda \subseteq \mu}\ c_{\lambda}\ \KL_{\lambda},
\]
where the coefficients $c_{\lambda}$ are monomials in $y_1,\ldots,y_{n-1}$ of degree at most 1 in each variable, and each KL element appears in the sum. These coefficients are recursively obtained by decomposing diagrams into ribbons, using equations \eqref{M2expansion} and \eqref{KLrecursion} below.
\end{theorem}

\begin{remark}
Note that since $\mu_k=k$ we have in this case that $u_k=v_k$ and $y_k=-\frac{[v_k]}{[v_k+1]}=-\frac{[u_k]}{[u_k+1]}$.
\end{remark}

\bigskip\noindent
Before proving this theorem, we first give an example.

\begin{example}[$n=3$]
Let 
\[
\M_{(2,1)}(u_1,u_2)
%\ =
%\begin{squarecells}{2}{2.2em}
%$\scriptstyle u_2+2$ & $\scriptstyle u_2+1$ \nl
%$\scriptstyle u_1+1$ \cr
%\cline{1-1}
%\end{squarecells}\ 
= T_2(u_1+1) T_4(u_2+1) T_3(u_2+2)\ \Delta_t(z_1,z_2,z_3) \Delta_t(z_4,z_5,z_6).
\]
The integral representation for $\M_{(2,1)}(u_1,u_2)$ is
\begin{align*}
\M_{(2,1)}(u_1,u_2) &=  \Delta_t(z_1,\ldots,z_6) \oint\oint \frac{1}{(2\pi\i)^{3 }} \Delta(w_1,w_2,w_3)\Delta_t (w_1,w_2,w_3) \phi_1(w_1)  \times\\
& \hphantom{\Delta_t(z_1,\ldots,z_N)} \left(1+y_1 S_{3}\right)\ \left(1+y_2 S_{5}\right)\cdot \phi_{3}(w_{2})\phi_{5}(w_{3}) \hspace{1mm}\d w_1 \d w_2  \d w_3,
\end{align*}
from which it immediately follows that
\begin{align*}
\M_{(2,1)}(u_1,u_2) &= \M_{(2,1)}(0,0) + y_1 \M_{(2)}(\cdot,0) + y_2 \M_{(1,1)}(0,0) + y_1 y_2 \M_{(1)}(\cdot,0),\\
&= \begin{squarecells}{2}{1em}
$\scriptstyle 2$ & $\scriptstyle 1$ \nl
$\scriptstyle 1$ \cr
\cline{1-1}
\end{squarecells} + 
y_1
\begin{squarecells}{2}{1em}
$\scriptstyle 2$ & $\scriptstyle 1$ \nl
\end{squarecells}
+ y_2
\begin{squarecells}{1}{1em}
$\scriptstyle 1$ \nl
$\scriptstyle 1$ \nl
\end{squarecells}
 + y_1y_2
\begin{squarecells}{1}{1em}
$\scriptstyle 1$ \nl
\end{squarecells}\\
&= \KL_{(2,1)} +y_1\KL_{(2)} +y_2 (\KL_{(1,1)}+\KL_{\emptyset}) +y_1y_2\KL_{(1)}.
\end{align*}
\end{example}
\bigskip

\begin{proof}[Proof of Theorem~\ref{th:expansion}]
The fact that the coefficients are polynomials in $y_1,\ldots,y_{n-1}$ of degree at most 1 in each variable follows immediately from the first equality in Theorem~\ref{DefMacdIntegral}. We will see below that they are in fact monomials.

We will use a recursive argument to prove that each KL element appears in the expansion. First note that the expansion of $\M_{n-1,n-2,\ldots,2,0}$ can be obtained from that of $M_{n-2,\ldots,1,0}$ by a shift of index $\KL_\lambda\rightarrow \KL_{\lambda+1^{n-2}}$, i.e. if
\be
\M_{n-2\ldots10} = \sum_{\lambda\subseteq \mu} c_{\lambda} \KL_{\lambda},
\label{M1expansion}
\ee
then
\be
\M_{n-1\ldots20} = \sum_{\lambda\subseteq \mu} c_{\lambda} \KL_{\lambda+1^{n-2}}.
\label{M2expansion}
\ee
This follows immediately by comparing the first line of Theorem~\ref{DefMacdIntegral} with and without shifted indices $\lambda_k \rightarrow \lambda_k+1$, $k=1,\ldots,n-1$.

Secondly, the polynomial $\M_{n-1\ldots1}$ is the image of $\M_{n-1\ldots20}$ under the action of $T_2(u_1+1)$, i.e.
\[
\M_{n-1\ldots1}(u_1,\ldots,u_{n-1}) = T_2(u_1+1)\M_{n-1\ldots20}(u_2,\ldots,u_{n-1}),
\]
and we can use the known action of $T_2(u_1+1)=T_2(1) + y_1$ on the KL polynomials appearing in the expansion of $\M_{n-1\ldots20}$. In order to explain this in terms of partitions, let us first define the notion of a Dyck ribbon.

Given a partition $\lambda$ and a diagonal touching an outer corner, a Dyck ribbon is a a strip $r_{\rm d}$ of boxes on the border of $\lambda$  from one outer corner to another, starting and ending on the chosen diagonal but never crossing it. The example in Figure~\ref{dyckribbon} illustrates this notion pictorially, see also \cite{ShigePZJ10}.

\begin{figure}[h]
\label{dyckribbon}
\centerline{\includegraphics[height=70pt]{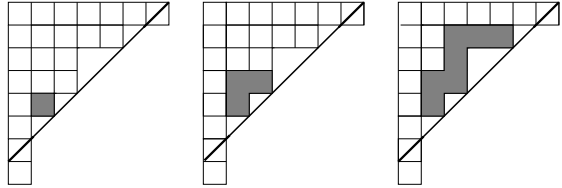}}
\caption{The partition $(7,5,3,3,2,1,1,1)$ with diagonal as shown has three Dyck ribbons.}
\end{figure}

Let $\lambda=(\lambda_1,\ldots,\lambda_k)$ be a partition with $k\leq n-2$, $\lambda'$ its dual, and let $\rho$ be the staircase $(\lambda_1'+1,\lambda_1',\ldots,1)$. Let $R_\rho$ denote the set of Dyck ribbons of $\lambda$ with diagonal the diagonal of $\rho$ that runs through its corner boxes. Then the action of $T_2(u_1+1)=T_2(1)+y_1$ on $\KL_{\lambda}$ is given by \cite{GierP07}

\be
\renewcommand{\arraystretch}{1.6}
T_2(u_1+1)\ \KL_{\lambda} = \left\{
\begin{array}{ll}
\left([2]+y_1\right) \KL_{\lambda} = \displaystyle\frac{[u_1+2]}{[u_1+1]} \KL_{\lambda} \qquad & k < n-2  \\
\KL_{\lambda,1} + y_1 \KL_{\lambda} + \sum_{r_{\rm d} \in R_\rho} \KL_{\lambda -r_{\rm d}} & k=n-2, 
\end{array}\right.
\label{KLrecursion}
\ee
where the last sum only appears if $\lambda_{n-2}=1$, and $\lambda-r_{\rm d}$ means the partition obtained by removing $r_{\rm d}$ from $\lambda$.

\begin{example}
For $n=4$,
\begin{align*}
T_2(u_1+1) \KL_{(2,2)} &= \KL_{(2,2,1)} + y_1 \KL_{(2,2)},\\
T_2(u_1+1) \KL_{(2,1)} &= \KL_{(2,1,1)} + y_1 \KL_{(2,1)} + \KL_{(2)} + \KL_{\emptyset}.\\
\end{align*}
\end{example}
\noindent
Equations \eqref{M2expansion} and \eqref{KLrecursion} give rise to a recursive procedure to determine the coefficients $c_{\lambda}$ in the expansion \eqref{M1expansion}. A hint of how the recursion works is obtained by considering the explicit the expansion of $\M_{(3,2,1)}$ using \eqref{KLrecursion}, given in \ref{ap:expansion}. 
\\

\noindent\textit{A recursion for $c_\lambda$}.\\
To complete the proof of Theorem~\ref{th:expansion},, we deduce a recursion for the coefficient $c_\lambda$ from \eqref{M2expansion} and \eqref{KLrecursion}, which is best explained in an example. We therefore show in \ref{backex} below how to recursively find the coefficient $c_{(2,2)}$ of $\KL_{(2,2)}$ in the expansion of $\M_{(6,5,4,3,2,1)}$. From this example it will be clear how to find the coefficient $c_\lambda$ for an arbitrary partition $\lambda$ in the expansion of any deformed maximal M polyomial.

The entire computation can be visualised in a single step, see Figure \ref{exsc22_overall}. In each iteration of the recursion, a maximal Dyck ribbon is attempted to being added to the smaller, inner, partition. This ribbon is not allowed to cross the diagonal formed by the diagonal of the staircase.  In the first iteration, we fill the first maximal\  Dyck ribbon with the numeral $1$ as the label of the iteration.  We then consider the smaller staircase obtained by ignoring the first column. In the second iteration, the second Dyck ribbon is added, and labelled with the numeral 2. This procedure continues for the six iterations. Where no Dyck ribbon can be added, a factor of $y_i$ is written in the diagram. Finally, the coefficient $c_{(2,2)}$ is equal to the product of these factors, $c_{(2,2)}=y_3y_5y_6$.

\begin{figure}[here!]
\begin{center}
\resizebox{5cm}{!}{\includegraphics{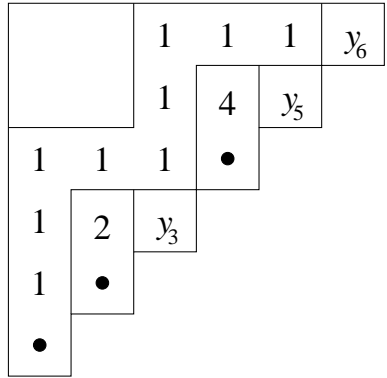}}
\caption{Summary of the coefficient computation
\label{exsc22_overall}}
\end{center}
\end{figure}

The reasoning above shows, as is also clear from the explicit example in ~\ref{backex}, that in this way one may compute the coefficient of $\KL_\lambda$ for any $\lambda$ in the expansion of $M_\rho$, with $\rho$ a staircase. The corresponding coefficient $c_\lambda$ is equal to the product of the factors $y_i$ on the diagonal, showing that this coefficient is indeed a monomial. It is also clear that each $c_\lambda$ is nonzero.

\end{proof}

\begin{cor}
Equation \eqref{M2KL} holds.
\end{cor}
\begin{proof}
First note that the case of equation \eqref{M2KL} corresponds to $y_i=\tau^{-1}=-[2]^{-1}$ for $i=1,\ldots,n$. Secondly, assigning signs as in Figure~\ref{cdef}, the signs within each Dyck ribbon in the recursion above add up to zero. The sum over signs between inner shape and staircase is therefore the same as the number of boxes containing a factor $y_i$, and hence the total coefficient $c_{\lambda}$ for each inner shape $\lambda$ is equal to $c_\lambda=\tau^{-n_\lambda}$, where $n_\lambda$ is the total sum of signed boxes.
\end{proof}

\section{Evaluations, constant term and Schubert polynomials}
\subsection{Evaluations}\label{eval}
The evaluations $z_i = 1$ of certain degenerate Macdonald and KL polynomials, normalised by dividing by the $t$-Vandermonde determinants, correspond to the number of (punctured) totally symmetric self complementary plane partitions (PTSSCPPs) with a weight in $\tau = -t -1/t$. Such enumerations were considered in \cite{DF06,DFZJ07,GierPZJ,FonsecaZJ08,FonsecaN10} using methods which we will discuss in Section~\ref{integrals}. Let us first give an example of two evaluations:

\begin{example}
\label{evaluations}
Let $\Delta  \Delta  =\Delta_t(z_1,z_2,z_3)\Delta_t(z_4,z_5,z_6)$. The ratio between the two extremal Macdonald polynomials satisfies
\be
\frac1{\Delta\Delta} M_{221100}\Big|_{z_i=1} = 3\tau^{-1}+3\tau+\tau^3.
\ee
A different result is obtained when specialising the sum of the KL basis, i.e. taking $c_\lambda=0$ in \eqref{M2KL}. 
\be
\frac1{\Delta\Delta} \left(KL_{111000}+KL_{110100}+KL_{101100}+KL_{110010}+KL_{101010}\right)\Big|_{z_i=1} = 1+3\tau+2\tau^2+\tau^3.
\ee
\end{example}

TSSCPPs, see e.g. \cite{MillsRR,Andrews94,Bressoud}, are in simple bijection to nonintersecting lattice paths consisting of north and north-east steps, and starting at positions $(i,-i)$ and ending on the $x$-axis. We further augment these paths, following \cite{DFZJ07}, with an extra step between the lines $y=0$ and $y=1$, in such a way that the difference of successive $x$ coordinates of the end positions are odd numbers, and the coordinate of the first endpoint is 1. We also assign a weight $\tau$ to vertical steps, and a weight $t$ to vertical augmented step.

\begin{example}\label{nilpaths}
The seven weighted paths for $n=2$ are
\bigskip

\centerline{\includegraphics[width=300pt]{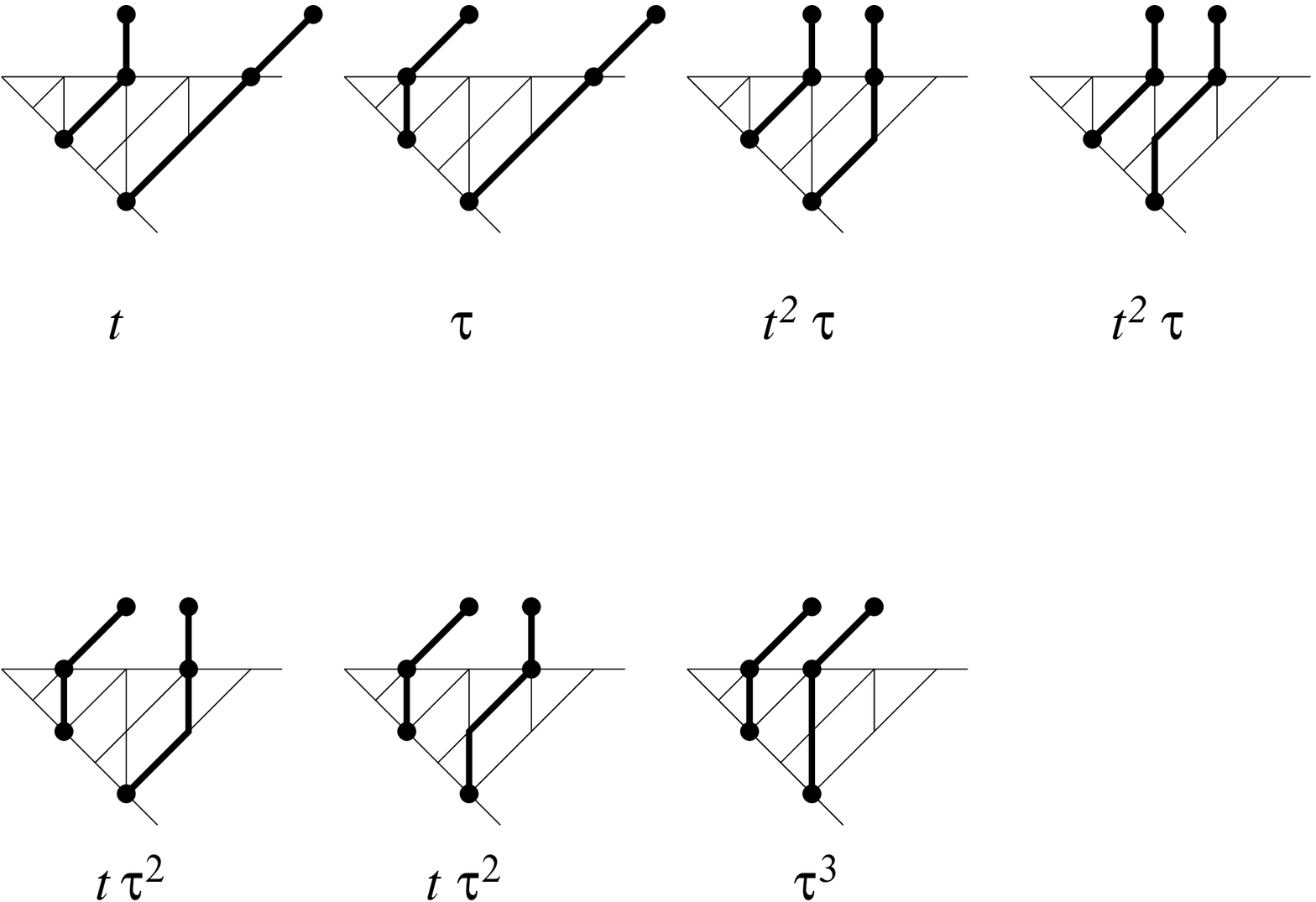}}

\bigskip\noindent
The weighted enumeration of these augmented paths is $N(t,\tau)=\tau^3 + 2t\tau^2+2t^2\tau+\tau+t$. 
\end{example}

Note that for special values of $t$ we obtain the evaluations in Example~\ref{evaluations}.
\[
N(\tau^{-1},\tau)=3\tau^{-1}+3\tau+\tau^3,\qquad N(1,\tau)=1+3\tau+2\tau^2+\tau^3.
\]
In the next section we clarify this correspondence and give a unified description of the M and KL polynomials.

\subsection{Constant term}\label{eval1}
In this section we further  investigate the homogeneous limit $z_i\rightarrow 1$. The following change of variables is useful
\be
z_i=\frac{1-t^{-1} \zeta_i}{1-t \zeta_i},\qquad w_i=\frac{1-t^{-1} x_i}{1-t x_i}.
\label{eq:changeofvar}
\ee
We also define
\[
\tau=-(t+t^{-1}).
\]
\begin{prop}[Homogenous limit]
\label{th:homlimit}
In the homogeneous limit $\zeta_i\rightarrow 0$ we obtain

\begin{multline*}
\M_{\lambda_{n-1},\ldots,\lambda_1}(u_1,\ldots,u_{n-1}) = \\ (t-t^{-1})^{n(n-1)}
\oint \cdots \oint \prod_{m=0}^{n-1}\frac{\d x_{m+1}\ (1+ y_m x_{m+1})}{2\pi\i\ x_{m+1}^{a_m}}\prod_{1\leq i<j\leq n} (x_j-x_i)(1+x_i x_j+\tau x_j),
\end{multline*}

where we recall that $a_k=\lambda_k+k+1$ and $y_k=-\frac{[v_k]}{[v_k+1]}, \ v_k=u_k+k-\lambda_k$ and $y_0=0$.
\end{prop}

\noindent 
A proof or Proposition~\ref{th:homlimit} is given in \ref{ap:homlimit}.

The contour integrals in Proposition~\ref{th:homlimit} pick out the residues around zero of each variable $x_i$. The integral can be rewritten as the following constant term expression
\be
A_{a_1,\ldots,a_{n-1}}(\tau,y_1,\ldots,y_{n-1}) = \CT \left( \prod_{i=0}^{n-1} (1 + y_ix_{i+1}) x_{i+1}^{-a_i+1}\prod_{i<j} (x_j - x_i)(1 + x_i x_j +\tau x_j) \right),
\label{CT1}
\ee
where $\CT =$ constant term with respect to each variable $x_i$. In the case $a_i=2i+1$ for $i=0,\ldots,n-1$\ and $y_1=\ldots = y_{n-1}=y$, expression \eqref{CT1} is equal to a constant term expression given by Di Francesco and Zinn-Justin in \cite{DFZJ07}. In the same paper, a formula (Formula 2.7) is given for the number of TSSCPPs according to the heights of the vertical steps. This formula reads
\be
N(t_0,\dots, t_{n-1}) = \CT \left( \prod_{i<j} \frac{(x_j - x_i)(1 + t_ix_j )}{1 - x_ix_j} \prod_i \frac{1 + t_0 x_i}{1 - x_i^2} \prod_{i=1}^n x_i^{-2i+2} \right)
\label{CT2}
\ee
When $t_1=\ldots t_{n-1}=\tau$, the two constant terms \eqref{CT1} and \eqref{CT2} are equal:

\begin{theorem}[Zeilberger-Di Francesco-Zinn Justin] 
\[
N(t,\tau,\ldots,\tau) = A_{1,3,\ldots,2n-1}(\tau,t,\ldots,t).
\]
\end{theorem}
\begin{proof}
This can be shown by symmetrising the argument of the constant term \eqref{CT1} with respect to the variables $x_i$, \cite{DFZJ07,Zeilberger}.
\end{proof}

We have not been able to establish a relation between $N$ and $A$ for general values of $y_i$ and $t_i$.
Fonseca and Zinn-Justin\cite{FonsecaZJ08}, and Fonseca and Nadeau \cite{FonsecaN10} consider the number of punctured totally symmetric self complementary plane partition with a puncture of size $r$. Their generation function is equal to
\[
N_{n,r}(\tau) = \CT \left( \prod_{i<j} \frac{(x_j - x_i)(1 + \tau x_j )}{1 - x_ix_j} \prod_i \frac{(1 + x_i)(1+\tau x_i)^{r}}{1 - x_i^2} \prod_{i=1}^n x_i^{-2i+2} \right)
\]

\begin{theorem}
Let $\mu^{(r)}=(\mu^{(r)}_{n-1},\ldots,\mu^{(r)}_1)$ be the staircase partition $(n-r-1,n-r-2,1,0,\ldots,0)$. Then we have

\[
N_{n-r,r}(\tau) = M_{\mu^{(r)}}(u_i=1-r;z_i=1) = A_{123\cdots,r-1,r+1,r+3,\ldots 2n-r-1}(\tau,1,\ldots,1).
\]

\end{theorem}

\begin{example}
For $n=5$ and $r=2$ we have $\mu_2=(2,1,0,0)$ and
\[
N_{3,2}(\tau)=M_{(2,1,0,0)}(u_1=\cdots=u_{4}=-1)=A_{12357}(\tau,1,1,1,1,1)=1+7\tau+12\tau^2+14\tau^3.
\]
\end{example}

\subsection{Schubert polynomials}

\emph{Schubert polynomials}  are a linear basis of the ring of polynomials.
  Though having been introduced in relation with Schubert calculus
on flag varieties, they will be specially appropriate for the computation of the functions
$N(t_0,\dots, t_{n-1})$.  Let us recall some facts about  these polynomials
(see e.g. \cite{Las03}).

Let $n$ be an integer, $\x=\{x_1,\ldots, x_n\}$,
$\y=\{y_1,y_2,\ldots,y_\infty\}$.
The \emph{Schubert polynomials} $Y_v(\x,\y)$, $v\in\mathbb{N}^n$,
are a linear basis of the space of polynomials in $\x$ 
with coefficients in $\y$.  They are recursively defined
starting from the case where $v$ is \emph{dominant}, i.e.
$v=\lambda=\lambda_1\ldots\lambda_n$, 
$\lambda_1\geq \lambda_2\ldots \geq \lambda_n\geq 0$
(we also say that $\lambda$ is a partition).
In that case
\[ Y_\lambda(\x,\y) := \prod_{i=1}^{n} \prod_{j=1}^{\lambda_i}  (x_i-y_j)   \, .\]
The general polynomials are then defined by
\be  
\partial_i Y_v(\x,\y) =
Y_{[\ldots, v_{i-1}, v_{i+1},v_i-1,v_{i+2},\ldots]}(\x,\y)\ ,
 \ v_i> v_{i+1}\, , 
\label{Schubertrecursion}
\ee
(in the case where $v_i \leq  v_{i+1}$, then $ \partial_i Y_v(\x,\y) =0$
since $\partial_i\partial_i=0$).

One can also index Schubert polynomials by permutations , the 
\emph{code of a permutation} \cite{Las03,Man01}) furnishing the correspondence between
the two indexings.

Since the image of $\partial_i$ is the space of polynomials
which are  symmetric in $x_i,x_{i+1}$,
the Schubert polynomial $Y_v(\x,\y)$ has this symmetry whenever
 $v_i \leq  v_{i+1}$. In particular, if $v$ is such that
 $v_1\leq \cdots\leq v_n$, then
$Y_v(\x,\y)$ is symmetric in  $x_1,\ldots,x_n$ and is called
a \emph{Gra{\ss}mannian  Schubert polynomial}.
It specializes to the Schur function in $\x$ indexed by the partition
$v_n\ldots v_1$ for $\y={\mathbf 0}=\{ 0,0,\ldots\}$. 
 
\begin{example}   A  
dominant Schubert polynomial can be visualised on the 
 diagram corresponding to the index of the polynomial.  Each box is filled with a factor of the form $(x_i-y_j)$, placed in the $i$th row and $j$th column.  The Schubert polynomial is equal to the product of these factors.  For $v=421$ the diagram is
\ben
\begin{squarecells}{4}{4em}
$x_1-y_1$&$x_1-y_2$&$x_1-y_3$&$x_1-y_4$ \nl
$x_2-y_1$& $x_2-y_2$ \cr \cline{1-2}
$x_3-y_1$ \cr\cline{1-1}
\end{squarecells}
\een 
so that we have
\ben
Y_{[4,2,1]}&=&(x_1-y_1)(x_1-y_2)(x_1-y_3)(x_1-y_4)(x_2-y_1)(x_2-y_2)(x_3-y_1).
\een
\end{example}

Following these preliminaries on Schubert polynomials, we now state
\begin{theorem}
\label{CT2Schu}
 Let $n$ be an integer, $\rho=n-1,\ldots,1$. 
For any $\lambda \subseteq \rho$, 
let $\mu=\lambda'$ be the conjugate partition, and let 
$(\rho/\lambda)=0,n-1-\mu_1,0,n-2-\mu_2,0,n-3-\mu_3,0,\dots$.

Then the constant term \eqref{CT2} is equal to
\bea
N(t_0,\dots,t_{n-1})=\sum_{\lambda \subseteq \rho,\hspace{1mm} \lambda \mbox{ {\small \it even}}}Y_{(\rho/\lambda)}({\bf \bar{y}},0),\label{eq:sfin}
\eea
and the alphabet ${\bf \bar{y}}$ being specialized into $\{t_0, t_1, 0, t_2, 0, t_3, 0,\dots\}$.
\end{theorem}

\begin{proof}[Proof of Theorem~\ref{CT2Schu}]
Let us analyse the expression \eqref{CT2}. One first extracts the kernel $\prod_{i<j} (1-x_i/x_j)$,
then the symmetric function $\prod_{i\le j} (1-x_ix_j )^{-1}$. What remains can be identified
with a Schubert polynomial. This allows us to use the following fundamental
scalar product on polynomials: given two polynomials, $f$, $g$, one defines
\bea\label{scalar}
\langle f,g \rangle=\CT \left(f(x_1,\dots,x_n)g(x_n^{-1},\dots,x_1^{-1})\prod_{1\le i <j \le n}(1-\frac{x_i}{x_j})\right).
\eea
It is proved in \cite{MFL09} that this scalar product  coincides with the usual scalar product when restricted to symmetric functions. The function $N(t_0,\dots, t_{n-1})$ is now written
\bea\label{nct}
N(t_0,\dots,t_{n-1}) &=&  \CT \left( \prod_{i<j} (1-\frac{x_i}{x_j})\prod_{i\leq j}(1-x_i x_j)^{-1} \prod_{i<j}(1+t_ix_j)\prod_{i=1}^{n} (1+ t_0 x_i) \prod_{i=1}^n x_i^{-i+1} \right)\nn \\
&=&  \left\langle \prod_{i\le j}(1-x_ix_j)^{-1},\prod_{i<j}(1+\frac{t_i}{x_{n+1-j}})
\prod_{i=1}^n (1+\frac{t_0}{x_i})x^{\rho_i}   \right\rangle \nn \\
&=&  \left\langle\prod_{i\le j}(1-x_ix_j)^{-1},
\prod_{i+j\le n}(x_i+t_j)
\prod_{i=1}^{n-1}(1+\frac{t_0}{x_i})\right\rangle
\eea
as the term $t_0/x_n$ does not contribute to the constant term. 

 Define the operator $\partial_\omega$ by taking a reduced decomposition $s_i\dots s_k$ of the permutation $(n,\dots,1)$ and putting  $ \partial_i\dots \partial_k =\partial_\omega $. 
Recall that the operator $\pi_{\omega} :\ f\mapsto\partial_\omega ( x^\rho f)$ is self adjoint and sends a
dominant monomial $x^{\lambda}$ onto the Schur function of the same index. Since the LHS of the scalar product is symmetric, one can symmetrize the RHS under $\pi_\omega$, transforming it into 

\bea
\partial_w \left(\prod_{i+j\le n}(x_i+t_j)
\prod_{i=1}^{n-1}(1+\frac{t_0}{x_i})x_i^{\rho_i}\right).
\eea

However
\ben
\prod_{i+j\le n}(x_i+t_j)\prod_i^{n-1}(1+\frac{t_0}{x_i})x_i^{\rho_i}
\een
is equal to the dominant Schubert polynomial of index $[2n-2,2n-4,\dots,2,0]$ for the alphabets ${\bf x}
=\{x_1,\dots,x_n\}$, ${\bf y}=\{-t_0,-t_1,0,-t_2,0,-t_3,\dots\}$. For example, for $n=4$, this Schubert polynomial is
\ben
\begin{squarecells}{6}{4em}
$x_1+t_0$ &$x_1 +t_1$& $x_1 +0$  & $x_1+t_2$ & $x_1 +0$ & $x_1 +t_3$ \nl
$x_2+t_0$ & $x_2+t_1$ & $x_2 +0$  & $x_2+t_2$ \cr \cline{1-4}
$x_3+t_0$ & $x_3+t_1$ \cr\cline{1-2}
\end{squarecells}
\een
We can therefore write
\bea
N(t_0,\dots,t_{n-1}) &=& \left\langle\prod_{i\le j}(1-x_ix_j)^{-1},\partial_\omega Y_{[2n-2,2n-4,\dots,2,0]}({\bf x},{\bf y})\right\rangle.
\eea
Using $\partial_\omega=\partial_{n-1}\circ\partial_{n-2}\partial_{n-1}\circ\dots \circ\partial_2\partial_3\dots \partial_{n-2}\partial_{n-1}$ and \eqref{Schubertrecursion}, we can simplify the RHS of the scalar product with
\ben
\partial_\omega Y_{[2n-2,2n-4,\dots,2,0]}({\bf x},{\bf y})
 = Y_{[0,1,\ldots, n-1]}({\bf x},{\bf y}).
\een

There is a Cauchy formula for Schubert polynomials \cite[(Theorem 10.2.6)]{Las03}, which decomposes a
Schubert polynomial into an alternating sum of products of Schubert polynomials
in each of the alphabets ${\bf x}$, ${\bf y}$ separately. In our case, $Y_{[0,1,...,n-1]}({\bf x},{\bf y})$ is
symmetric in $x_1,\dots, x_n$ and decomposes into a sum of Schur functions in ${\bf x}$:
\bea
Y_{[0,1,2,\dots,n-1]}({\bf x},{\bf y})=\sum_{\lambda \subseteq \rho}S_\lambda({\bf x})Y_{(\rho/\lambda)}({\bf \bar{y}},0),
\eea
the definition of the index $(\rho/\lambda)$ and the alphabet ${\bf \bar{y}}$ being given in the statement of the theorem.

The LHS of the scalar product is equal to the sum of all Schur functions
indexed by an even partition (i.e. with even parts). Therefore, the scalar product is equal to
\ben
N(t_0,\dots,t_{n-1})=\sum_{\lambda \subseteq \rho,\hspace{1mm} \lambda \mbox{ {\small even}}}Y_{(\rho/\lambda)}({\bf \bar{y}},0).  \\
\een
\qedhere
\end{proof}

\begin{example}
 \item For $n=3$, we have $\rho=21$. 
 The two possible even sub-partitions are $\lambda=\emptyset$ and $\lambda=2$.  The two conjugate partitions are $\mu=\lambda'=\emptyset$, and $\mu=\lambda'=11$ respectively.  Considering the first of these we have
 \ben
 (\rho/\lambda)=(21/0)=[0,n-1-\mu_1,0,n-2-\mu_2,0]=[0,2,0,1,0],
 \een
 and for the second sub-partition we have
 \ben
 (\rho/\lambda)=(21/2)=[0,n-1-\mu_1,0,n-2-\mu_2,0]=[0,1,0,0,0].
 \een
 Thus the two polynomials contributing to the constant term $N(t_0,t_1,t_2)$ are
 \ben
 Y_{[0,2,0,1,0]}({\bf \bar{y}},0) \mbox{ and }  Y_{[0,1,0,0,0]}({\bf \bar{y}},0) 
 \een
\end{example}
%
%%%%%%%%%%%%%%%%%%%%%%%%%%%%%
\subsection{Determinant expression}
%%%%%%%%%%%%%%%%%%%%%%%%%%%%%
The Schubert polynomials in Theorem~\ref{CT2Schu} have a determinantal expression that we are now going to introduce. The complete symmetric function of degree $k$ is defined to be the sum of all monomials of total degree $k$ -- see e.g. \cite{Mac95}. Let $h_k(r)$ denote the complete symmetric function of degree $k$ in the alphabet ${\bf \bar{y}}=\{y_1,\dots,y_{2r-2}\}=\{t_0,t_1,0,t_2,0,t_3,\dots,t_{r-1}\}$.  We also use the notation $\mu=\lambda'$. 

The indices of the Schubert polynomials in ${\bf \bar{y}}$ on the RHS of \eqref{eq:sfin}, are all codes which can be obtained from an increasing partition under the Schubert recursion \eqref{Schubertrecursion}. Therefore, these polynomials are all images of Schur functions in ${\bf \bar{y}}$ under divided differences.

Schur functions have a determinantal expression in terms of complete functions (Jacobi-Trudi determinants).
We need to generalize these determinants by allowing flags of alphabets.
Given two partitions $\rho,\lambda$, 
an alphabet ${\bf \bar{y}}=\{y_1,y_2,\dots\}$,
an increasing sequence of positive integers $\Phi=(\phi_1,\phi_2,\dots)$, 
then the flag  Schur function $S_{\rho/\lambda}({\mathbf{ \bar y}},\Phi)$
is equal to 
\[
S_{\rho/\lambda}({\mathbf{ \bar y}},\Phi)= \det\left(h_{\rho_i-\mu_j-i+j}({\phi_i})\right)_{i,j=1\dots n-1}.
\]
The action of divided differences on such determinants, under some hypotheses
which are satisfied in our case, is easy, see e.g. \cite[Lemma 1.4.5]{Las03} or \cite[Corollary 2.6.10]{Man01}. At each step any divided difference acts on a single row only (or on no row). From \eqref{Schubertrecursion} it follows that this action consists in decreasing the indices of the complete symmetric functions in this row and transforming their argument $r$ into $r+1$. We thus obtain

\begin{prop}
Let $\mu=\lambda'$ be the partition conjugate to $\lambda$. Then
\bea
Y_{(\rho/\lambda)}({\bf \bar{y}},0) = S_{\rho/\lambda}({\mathbf{ \bar y}},\Phi) 
\eea
where $\Phi=(2,4,\dots, 2n-2)$. Equivalently, suppressing the 0's in ${\bf \bar{y}}$ and using the alphabet ${\bf \tilde{y}}=(t_0,t_1,t_2,t_3,\ldots)$ with $\tilde{h}_i(r)$ the complete symmetric function in the alphabet ${\bf \tilde{y}}$, one has
\bea
Y_{(\rho/\lambda)}({\bf \tilde{y}},0) = S_{\rho/\lambda}({\bf \tilde{y}},{\widetilde{\Phi}})=\det\left(\tilde{h}_{\rho_i-\mu_j-i+j}({\tilde{\phi}_i})\right)_{i,j=1\dots n-1}
\label{detexp}
\eea
where the flag is equal to $\widetilde{\Phi}=(2,3,\ldots,n)$. 
\end{prop}

We illustrate this is in the following examples.
\begin{example}
\[  
Y_{[4,6,7]}({\bf \bar{y}},0) = 
\begin{vmatrix}   
h_4(1) & h_7(1) & h_9(1) \\
h_3(2) & h_6(2) & h_8(2) \\
h_2(3) & h_5(3) & h_7(3)  
\end{vmatrix}. 
\]
We transform this into $Y_{[0,3,0,4,0,4]}$ using the sequence of divided difference operators $\partial_1\partial_3\partial_2\partial_5\partial_4\partial_3$. The first step gives
\[ Y_{[4,6,7]} \stackrel{\partial_3}{\to} Y_{[4,6,0,6]} = 
\begin{vmatrix}   
h_4(1) & h_7(1) & h_9(1) \\
h_3(2) & h_6(2) & h_8(2) \\
h_1(4) & h_4(4) & h_6(4)   
\end{vmatrix},
\]
and the remaining transformations lead to \\
\[
Y_{[4,6,0,6]}\stackrel{\partial_4}{\to}  Y_{[4,6,0,0,5]} \stackrel{\partial_5}{\to}  Y_{[4,6,0,0,0,4]} \stackrel{\partial_2}{\to}  Y_{[4,0,5,0,0,4]} \stackrel{\partial_3}{\to}  Y_{[4,0,0,4,0,4]}  \stackrel{\partial_1}{\to}  Y_{[0,3,0,4,0,4]} =
\begin{vmatrix}   
h_3(2) & h_6(2) & h_8(2) \\
h_1(4) & h_4(4) & h_6(4) \\
0 & h_2(6) & h_4(6) \end{vmatrix}. 
\]
Recall that we have taken the alphabet ${\bf \bar{y}}=(t_0,t_1,0,t_2,\ldots)$, and the flag consists of exactly the positions of the non-zero components of the index. Then
\[
Y_{[0,3,0,4,0,4]} = 
\begin{vmatrix}   
\tilde{h}_3(2) & \tilde{h}_6(2) & \tilde{h}_8(2) \\
\tilde{h}_1(3) & \tilde{h}_4(3) & \tilde{h}_6(3) \\
0 & \tilde{h}_2(4) & \tilde{h}_4(4) \end{vmatrix}. 
\]
\end{example}

\begin{example}\label{detex}
Using \eqref{detexp}, we write
\ben
Y_{[0,2,0,1,0]}({\bf \bar{y}},0) &=& \dt{h_2(2)&h_3(2)\\h_0(3)&h_1(3)} \\
&=& \dt{(t_0^2+t_0t_1+t_1^2)&(t_0^3+t_0^2t_1+t_0t_1^2+t_1^3)\\1&(t_1+t_1+t_2)} \\
&=& t_0^2t_1+t_0^2t_2+t_0t_1^2+t_0t_1t_2+t_1^2t_2.
\een
Similarly
\ben
Y_{[0,1,0,0,0]}({\bf \bar{y}},0) &=& \dt{h_1(2)&h_2(2)\\h_{-1}(3)&h_0(3)} \\
&=& \dt{(t_0+t_1)&(t_0^2+t_0t_1+t_1^2)\\0&1} \\
&=&t_0+t_1.
\een\\
and thus
\bea\label{exschu}
N(t_0,t_1,t_2)=Y_{[0,2,0,1,0]}({\bf \bar{y}},0)+Y_{[0,1,0,0,0]}({\bf \bar{y}},0)=t_0+t_1+t_0^2t_1+t_0^2t_2+t_0t_1^2+t_0t_1t_2+t_1^2t_2.
\eea
\end{example}

\subsection{Combinatorial interpretation}
Having interpreted the constant term in terms of Schubert polynomials, one gets for free a combinatorial
interpretation in terms of tableaux. Indeed, $Y_{[0,n-1,0,\dots,3,0,2,0,1]}({\bf \bar{y}})$ is equal to a sum of semi-standard tableaux of staircase shape, that we represent in the French way in the Cartesian plane, satisfying a flag condition.  Taking the alphabet ${\bf \tilde{y}}=(t_0,t_1,t_2,t_3,\dots)$, the Schubert polynomial becomes $Y_{[0,n-1,n-2,\dots,2,1]}({\bf \tilde{y}})$, which can be interpreted as the sum of tableaux of staircase shape such that the bottom row belongs to $\{ 1,2\}^*$, i.e only $t_0$ and $t_1$ can be used as fillings for this row, the next one to $\{ 1,2,3\}^*,\dots,$ the top one to $\{ 1,\dots,n\}^*$. We say that such a tableau satisfies the flag condition $2,3,\dots,n$.
For $n=4$, this is

\[
\begin{array}{ll}
\raisebox{-5mm}{
\begin{invsquarecells}{3}{2.5em} 
\cline{1-1} 
\cr \cline{1-2}
&\cr \cline{1-3}
&&  \\
\cline{1-3}
\end{invsquarecells}}
\qquad &
\begin{array}{l}
{\rm flag}\\[3.5mm]
4 \\[5mm]
3 \\[5mm]
2\hspace{5.5cm} 
\end{array}
\end{array}
\]

This interpretation remains valid for skew shapes.  As usual in the theory of tableaux, this result is obtained by introducing an extra alphabet of small letters, which fill the inner shape.  The valuation of these tableaux for a fixed inner shape with even columns are exactly the determinants  written above.  \\
Therefore, (\ref{eq:sfin}) becomes
\begin{theorem}
The constant term $N(t_0,\dots,t_{n-1})$ is the commutative image of the sum of tableaux of outer shape $\rho=n-1,\dots,1$, inner shape with columns of even lengths, satisfying the flag condition $2,3,\dots,n$.
\end{theorem}
In fact, one can reprove directly that the generating function $N(t_0,\dots,t_{n-1})$  is given by a sum of tableaux. \\ \\

\begin{example}
 For $n=3$, the sum 
\[
Y_{(21/0)}({\bf \bar{y}},0) + Y_{(21/2)}({\bf \bar{y}},0) =
   Y_{[0,2,0,1,0]}({\bf \bar{y}},0) + 
    Y_{[0,1,0,0,0]}({\bf \bar{y}},0)
\]
is given by the  tableaux 
\ben
\begin{invsquarecells}{2}{2.5em} 
\cline{1-1}
$t_2$\cr\cline{1-2}
$t_0$&$t_0$ \nl
\end{invsquarecells} \hspace{4mm}
\begin{invsquarecells}{2}{2.5em} 
\cline{1-1}
$t_2$\cr\cline{1-2}
$t_0$&$t_1$ \nl
\end{invsquarecells}  \hspace{4mm}
\begin{invsquarecells}{2}{2.5em} 
\cline{1-1}
$t_2$\cr\cline{1-2}
$t_1$&$t_1$ \nl
\end{invsquarecells}  \hspace{4mm}
\begin{invsquarecells}{2}{2.5em} 
\cline{1-1}
$t_1$\cr\cline{1-2}
$t_0$&$t_0$ \nl
\end{invsquarecells}\hspace{4mm}
\begin{invsquarecells}{2}{2.5em} 
\cline{1-1}
$t_1$\cr\cline{1-2}
$t_0$&$t_1$ \nl
\end{invsquarecells}\hspace{4mm}
\begin{invsquarecells}{1}{2.5em} 
 \cline{1-1}
$t_0$\cr\cline{1-1}
\end{invsquarecells}\hspace{4mm}
 \begin{invsquarecells}{1}{2.5em} 
 \cline{1-1}
$t_1$\cr\cline{1-1}
\end{invsquarecells}
\een \\
and this agrees with
\ben
N(t_0,t_1,t_2)=t_0^2t_2+t_0t_1t_2+t_1^2t_2+t_0^2t_1+t_0t_1^2+t_0+t_1.
\een 
\end{example}
\bigskip

We now introduce another object which is in bijection with TSSCPP.
We take the staircase shape now in the north-west corner (here for $n=5$):
\setlength{\columnsep}{1pt}
\begin{multicols}{3}{
\ben
\hspace{3.5cm}&& \Box \hspace{5mm} \Box \hspace{5mm}\Box \hspace{5mm} \Box \\
&& \Box \hspace{5mm} \Box \hspace{5mm}\Box  \\
&& \Box \hspace{5mm} \Box \\
&& \Box
\een
\\
\vspace{0.8cm}\begin{center} \hspace{1.9cm} with border:\vspace{1.5cm} \end{center}
\ben
\\[1mm]
&&\Box \hspace{5mm} \Box \hspace{5mm}\Box \hspace{5mm} \Box \hspace{4.4mm} (0)\hspace{2cm}\\
&&\Box \hspace{5mm} \Box \hspace{5mm}\Box  \hspace{4.4mm} (1)\\
&&\Box \hspace{5mm} \Box \hspace{4.4mm} (2)\\
&&\Box \hspace{4.4mm} (3)
\een
}
\end{multicols}
\noindent
The TSSCPP are in bijection with the fillings of the completed staircase such that rows are weakly decreasing from left to right, and decrease by $1$ at most (taking into account the border).
Each filling corresponds to a configuration of non-intersecting lattice paths, beginning at $(i,-i)$, obtained by reading each row from right to left, and an increase in the integer corresponds to a vertical step. \\ 
\begin{example} For the filling
\ben
\hspace{1cm}&& \hspace{10mm}1\hspace{5mm} 1\hspace{5mm}1 \hspace{5mm}0  \hspace{3.7mm} (0) \\
&& \hspace{10mm}2 \hspace{5mm}2 \hspace{5mm}1  \hspace{3.7mm} (1)  \\
&& \hspace{10mm}3 \hspace{5mm} 3  \hspace{3.7mm} (2) \\
&& \hspace{10mm}3 \hspace{3.7mm} (3) 
\een
we obtain the non-intersecting lattic paths (NILP) shown in Figure \ref{nilpex}.
\begin{figure}[here!]
\begin{center}
\resizebox{6cm}{!}{\includegraphics{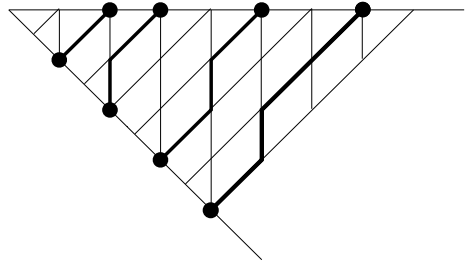}}
\caption{Non-intersecting lattice path configuration. Reading each row from right to left, an increase in value corresponds to a vertical step\label{nilpex}}
\end{center}
\end{figure}\\ 
\end{example}
Now add a left column (column $0$) as follows.  If the number of rows is even ({\it resp. odd}), then it consists of the even ({\it resp. odd}) numbers immediately bigger than or equal to the numbers in the first column. 
For example,
\[
\begin{array}{ccc}
\renewcommand{\arraystretch}{1.2}
\begin{array}{lllll}
\hspace{10mm}1\hspace{5mm} 1\hspace{5mm}1 \hspace{5mm}0  \hspace{3.7mm} (0) \\
\hspace{10mm}2 \hspace{5mm}2 \hspace{5mm}1  \hspace{3.7mm} (1)  \\
\hspace{10mm}3 \hspace{5mm} 3  \hspace{3.7mm} (2) \\
\hspace{10mm}3 \hspace{3.7mm} (3) 
\end{array}
&
{\rm becomes\;} 
&
\renewcommand{\arraystretch}{1.2}
\begin{array}{lllll}
2\hspace{5mm} 1\hspace{5mm} 1\hspace{5mm}1 \hspace{5mm}0  \hspace{3.7mm} (0) \\
2 \hspace{5mm} 2\hspace{5mm}2 \hspace{5mm}1  \hspace{3.7mm} (1)  \\
4 \hspace{5mm} 3\hspace{5mm} 3  \hspace{3.7mm} (2) \\
4 \hspace{5mm} 3\hspace{3.7mm} (3) 
\end{array}
\end{array}
\]
This corresponds to the completion defined by Di Francesco and Zinn-Justin in \cite{DFZJ07}, and discussed in Section \ref{eval} (see Example \ref{nilpaths}).  Again the rows describe the successive paths (read from right to left) composing the TSSCPP, an increase corresponding to a vertical step.\\
\begin{example} For the filling, with completion described above, we have the path configuration shown in Figure \ref{nilpcompex}.\\
\begin{figure}[here!]
\begin{center}
\resizebox{6cm}{!}{\includegraphics{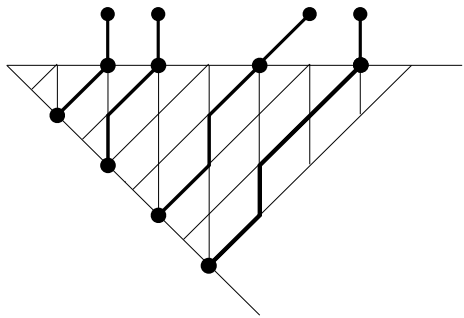}}
\caption{Augmented non-intersecting lattice path configuration. Again reading each row from right to left, and an increase in value corresponds to a vertical step\label{nilpcompex}}
\end{center}
\end{figure}
\end{example}
But such an object can also be read as a usual skew Young tableau (in the Cartesian plane, with strictly decreasing columns, and weakly increasing rows)  of outer shape the staircase, with inner shape the diagram with lengths of columns the complement of column $0$ to $(n-1)^{n-1} $. In the example above the inner shape would be $4444 - 2244 = 2200$. One records the positions of the decreases in row $1, 2,\dots,$ (reading left to right), and this 
gives columns $1, 2,\dots$ of the tableau (reading bottom to top). Such a tableau satisfies the flag condition 
$2, 3,\dots$, its commutative evaluation is exactly the weight of the corresponding 
TSSCPP (where $t_0$ is the weight of each augmented vertical step, $t_1$ the weight of each vertical step in the row immediately below this, etc.). The example becomes 
\ben
 \begin{invsquarecells}{4}{2.5em} 
 \cline{1-1} 4
 \cr \cline{1-2}
1&3\cr \cline{1-3}
$\cdot$&$\cdot$&3 \cr\cline{1-4}
$\cdot$&$\cdot$&1&1 \cr \cline{1-4}
\end{invsquarecells}\\
\een
and thus has a weight $t_0^3t_2^2t_3$. \\ \\
In summary, each weighted NILP corresponds to a filling, as described above.  Each filling corresponds to a skew tableau, satisfying a flag condition, describing a monomial of the Schubert polynomial $Y_{(\rho/\lambda)}({\bf \bar{y}},0)$ for some even partition $\lambda$.  
%%%%%%%%%%%%%%%%%%%%%%%%%%%%%%%%%%%%%%%%%%%%%%%%%
\section*{Acknowledgement}
AL and JdG would like to thank the hospitality of the Mathematisches Forschungsinstitut Oberwolfach, where part of this work was done, and we have greatly benefited from discussions with Tiago Fonseca, Pavel Pyatov and Keiichi Shigechi. We thank the Australian Research Council (ARC) for financial support.

\appendix

\section{Explicit example of the recursion in the proof of Theorem~\ref{th:expansion}}
\label{backex}
\begin{figure}[here!]
\begin{center}
\subfigure[The maximal Dyck ribbon is shaded]{\label{exsc22_1a}\resizebox{3cm}{!}{\includegraphics{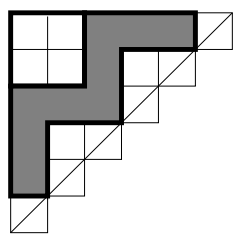}}}\qquad\qquad
\subfigure[The first column is removed, and new inner shape drawn]{\label{exsc22_1b}\resizebox{3cm}{!}{\includegraphics{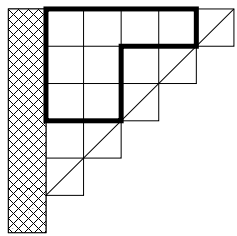}}}
\end{center}
\caption{Steps 1 and 2 in the first iteration}
\end{figure}

Consider the Macdonald polynomial $\M_{(6,5,4,3,2,1)}$, i.e. corresponding to the maximal staircase of size $n=7$.  We wish to find the coefficient of $\KL_{(2,2)}$ in the expansion of this polynomial. To do so, we begin by drawing the partition $(2,2)$ inside the staircase. \\
\begin{figure}[here!]
\begin{center}
\resizebox{3cm}{!}{\includegraphics{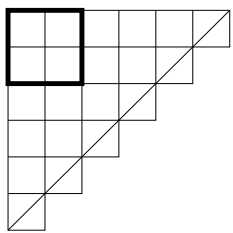}}
\caption{The maximal staircase of size $n=7$, with inner shape $(2,2)$
\label{exsc22}}
\end{center}
\end{figure}\\
We now show the recursion step by step. \\ \\
$\bullet$ {\rm First iteration}\medskip\\
{\rm Step 1}:
We draw in the diagonal bounding the staircase, and add the maximal Dyck ribbon to the inner shape (if possible).  This is shown in Figure~\ref{exsc22_1a}. \\
{\rm Step 2}:
The second step is to delete the first column.  Our inner shape is now the remainder of $(2,2)$ and the added Dyck ribbon, i.e. the new inner shape becomes $(4,2,2)$, see Figure~\ref{exsc22_1b}. 

\begin{figure}[here!]
\begin{center}
\subfigure[The maximal Dyck ribbon is just the shaded single box]{\label{exsc22_2a}\resizebox{2.5cm}{!}{\includegraphics{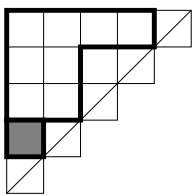}}}\qquad\qquad
\subfigure[The first column is again removed, and new inner shape drawn]{\label{exsc22_2b}\resizebox{2.5cm}{!}{\includegraphics{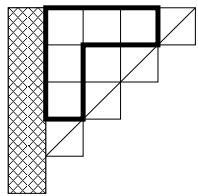}}}
\end{center}
\caption{Steps 1 and 2 in the second iteration}
\end{figure}

\bigskip\noindent
$\bullet$ {\rm Second iteration} \medskip\\
We now repeat steps 1 and 2, with the new inner shape and smaller staircase: \\
{\rm Step 1}:
We draw in the diagonal bounding the staircase, and add the maximal Dyck ribbon to $(4,4,2,2)$ inside $(5,4,3,2,1,1)$. This is shown in Figure \ref{exsc22_2a}. \\
{\rm Step 2}: We again remove the first column, see Figure \ref{exsc22_2b}, the new inner shape is $(3,1,1)$.

\begin{figure}[here!]
\begin{center}
\subfigure[No Dyck ribbon can be added]{\label{exsc22_3a}\resizebox{2cm}{!}{\includegraphics{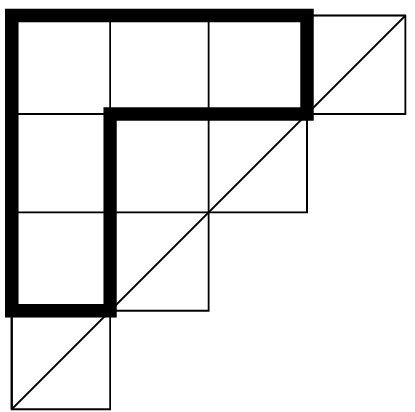}}}\qquad\qquad
\subfigure[The first column is again removed]{\label{exsc22_3b}\resizebox{2cm}{!}{\includegraphics{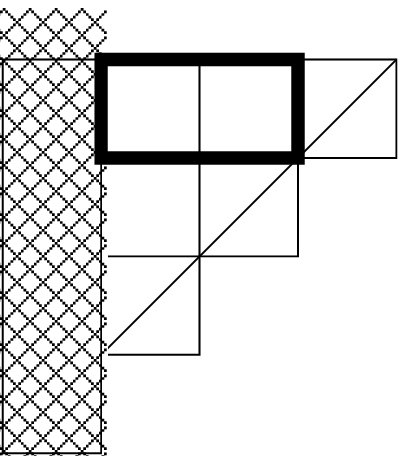}}}
\end{center}
\caption{Steps 1 and 2 in the third iteration}
\end{figure}

\bigskip\noindent
$\bullet$ {\rm Third iteration} \medskip\\
{\rm Step 1}:
It is not possible to add a Dyck ribbon, see Figure \ref{exsc22_3a}. This has the consequence that now we pick up a factor of $y_3$. \\
{\rm Step 2}:
We again remove the first column, see Figure \ref{exsc22_3b}. The new inner shape is $(2)$.

\begin{figure}[here!]
\begin{center}
\subfigure[Maximal Dyck ribbon]{\label{exsc22_4a}\resizebox{1.5cm}{!}{\includegraphics{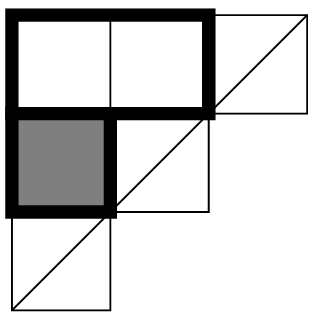}}}\qquad\qquad
\subfigure[First column removed]{\label{exsc22_4b}\resizebox{1.5cm}{!}{\includegraphics{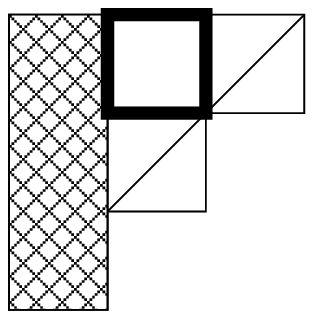}}}
\end{center}
\caption{Steps 1 and 2 in the fourth iteration}
\end{figure}

\bigskip\noindent
$\bullet$ {\rm Fourth iteration} \medskip\\
We repeat steps 1 and 2, with the new inner shape and smaller staircase: \\
{\rm Step 1}:
We draw in the diagonal bounding the staircase, and add the maximal Dyck ribbon, see Figure \ref{exsc22_4a}. \\
{\rm Step 2}:
We again remove the first column, see Figure \ref{exsc22_4b}. The new inner shape is $(1)$.

\begin{figure}[here!]
\begin{center}
\subfigure{\label{exsc22_5a}\resizebox{1cm}{!}{\includegraphics{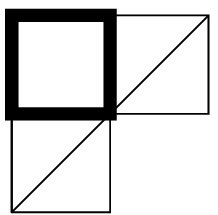}}}\qquad\qquad
\subfigure{\label{exsc22_5b}\resizebox{1cm}{!}{\includegraphics{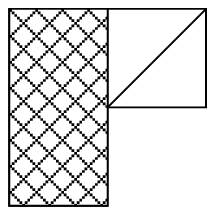}}}\qquad\qquad
\subfigure{\label{exsc22_6a}\resizebox{0.5cm}{!}{\includegraphics{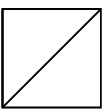}}}\qquad\qquad
\subfigure{\label{exsc22_6b}\resizebox{0.5cm}{!}{\includegraphics{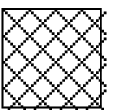}}}
\end{center}
\caption{Steps 1 and 2 in the fifth and sixth iteration\label{exsc22_5}}
\end{figure}

\bigskip\noindent
$\bullet$ {\rm Fifth and sixth iterations} \medskip\\
Repeating steps we see that no further Dyck ribbons can be added, see Figure~\ref{exsc22_5}. For these iterations therefore result in factors $y_5$ and $y_6$: \\

\noindent
Since we could not add Dyck ribbons at iterations 3, 5 and 6, we picked up corresponding factors $y_3$, $y_5$ and $y_6$. The claim is that the coefficient of $KL_{(2,2)}$ in the expansion of $M_{(6,5,4,3,2,1)}$ is just $y_3y_5y_6$.

\section{Proof of Theorem~\ref{DefMacdIntegral}}
\label{se:DefMacdIntegralProof}

\subsection{One row} 
First we recall
\be
\phi_i(w) = \prod_{m=1}^i \frac{1}{w-z_m} \prod_{m=i+1}^N\frac{1}{tw-t^{-1}z_m},
\ee
and we have

\begin{align}
T_i(v) \Delta_t(z_1,\ldots,z_N)\phi_{i}(w) &=  \Delta_t(z_1,\ldots,z_N) \Big(\phi_{i-1}(w) + \phi_{i+1}(w) + \frac{[1-v]}{[v]} \phi_i(w)\Big),
\label{eqmotiona}\\
T_i(v) \Delta_t(z_1,\ldots,z_N)\phi_{j}(w) &= \frac{[1+v]}{[v]} \Delta_t(z_1,\ldots,z_N)\phi_{j}(w)\qquad (j\neq i).
\label{eqmotionb}
\end{align}
We now look for a solution $ \psi _i$ of \eqref{eqmotiona} and \eqref{eqmotionb} satisfying the boundary conditions $ \psi _0= \psi _{N+1}=0$. These boundary conditions can be fixed using Cauchy's theorem. Let $C_0$ denote the contour lying around all poles at $z_i$ ($i=1,\ldots,N$), and let $C_1$ be the contour lying around the poles at $t^{-2}z_i$. As $\phi_i(w)\d w$ has no pole at infinity, the contour $C_0$ may be deformed to $C_1$, and we have
\be
\oint_{C_0} \frac{1}{2\pi\i}\ \phi_i(w) \hspace{1mm}\d w=-\oint_{C_1} \frac{1}{2\pi\i}\ \phi_i(w)\hspace{1mm}\d w.
\ee
By Cauchy's theorem, the first integral obviously gives zero for $i=0$, and second is clearly zero for $i=N$. 

Let us now define
\be
 \psi _{i} = \Delta_t(z_1,\ldots,z_N) \oint_{C_0} \frac{1}{2\pi\i}\ \phi_i(w)\hspace{1mm}\d w,
\ee
and note that 
\[
 \psi _2(v_0) := T_1(v_0+1)  \psi _1 =  \psi _2 - \frac{[v_0]}{[v_0+1]} \psi _1. 
\]
Let us further introduce the shift operator $S_i$ which acts as
\be
S_i \phi_i = \phi_{i-1},\qquad S_i \phi_j = \phi_j\quad (j\neq i).
\label{shiftoperator}
\ee
We then have that
\begin{align*}
 \psi _2(v_0) &= \Delta_t(z_1,\ldots,z_N)\oint  \frac{1}{2\pi\i}\ \left(1-\frac{[v_0]}{[v_0+1]} S_2\right)\ \phi_2(w)\hspace{1mm}\d w\\
&=\Delta_t(z_1,\ldots,z_N)\oint \frac{1}{2\pi\i}\ \left(\frac{1}{[v_0+1]} \frac{t^{v_0+1} w -t^{-v_0-1}z_2 }{tw-t^{-1}z_2} \right)\phi_2(w)\hspace{1mm}\d w.
\end{align*}
More generally we define
\be
 \psi_{a_0}(v_0)=T_{a_0-1}(v_0+1) \cdots T_1(v_0+a_0-1) \psi_1 ,
\label{psikdef}
\ee
which we denote by 
\[
 \psi_{a_0}(v_0)= 
\begin{squarecells}{6}{3.5em}
$\scriptstyle  v_0+a_0-1$ &  &  &  & $\scriptstyle v_0+2$ & $\scriptstyle v_0+1$ \nl
\end{squarecells}\ .
\]

\begin{lemma}
\label{lem:onerow}
For $a_0\geq 1$,
\begin{align*}
 \psi_{a_0}(v_0) &= \Delta_t(z_1,\ldots,z_N)\oint  \frac{1}{2\pi\i}\ \left(1-\frac{[v_0]}{[v_0+1]} S_{a_0}\right)\ \phi_{a_0}(w)\hspace{1mm}\d w\\
&=\Delta_t(z_1,\ldots,z_N)\oint \frac{1}{2\pi\i}\ \left(\frac{1}{[v_0+1]}\frac{t^{v_0+1} w -t^{-v_0-1}z_{a_0}}{tw-t^{-1}z_{a_0}} \right)\phi_{a_0}(w)\hspace{1mm}\d w.
\end{align*}
\end{lemma}
\begin{proof}
The proof is by induction. From \eqref{psikdef} it follows that
\[
 \psi_{a_0+1}(v_0) = T_{a_0}(v_0+1) \psi_{a_0}(v_0+1).
\] 
Using the properties \eqref{eqmotiona} and \eqref{eqmotionb} it is then easy to see that 

\begin{align*}
\psi_{a_0+1}(v_0) &= \oint  \frac{1}{2\pi\i}\ T_{a_0}(v_0+1) \Delta_t(z_1,\ldots,z_N)\left( \phi_{a_0}(w) - \frac{[v_0+1]}{[v_0+2]} \phi_{a_0-1}(w)\right)\hspace{1mm}\d w \\
&=  \Delta_t(z_1,\ldots,z_N)\oint  \frac{1}{2\pi\i}\ \left( \phi_{a_0-1}(w) + \phi_{a_0+1}(w) - \frac{[v_0]}{[v_0+1]} \phi_{a_0}(w) -\phi_{a_0-1}(w)\right)\hspace{1mm}\d w\\
&=  \Delta_t(z_1,\ldots,z_N)\oint  \frac{1}{2\pi\i}\ \left( 1 - \frac{[v_0]}{[v_0+1] }S_{a_0+1}\right) \phi_{a_0+1}(w)\hspace{1mm}\d w.
\end{align*}

\end{proof}
\bigskip

\subsection{Two rows}
Now we consider the function
\be
\phi_{a_0 a_1}(w_1,w_2) = \Delta_t(z_1,\ldots,z_N)\ A(w_1,w_2)\ \phi_{a_0}(w_1) \phi_{a_1}(w_2).
\label{tworowansatz}
\ee
We then find that 
\[
T_2( v_1 ) \phi_{12}(w_1,w_2) = \phi_{11}(w_1,w_2)+\phi_{13}(w_1,w_2)+ \frac{[1-v_1]}{[v_1]} \phi_{12}(w_1,w_2).
\]
The ``unwanted term'' $\phi_{11}(w_1,w_2)$ can be made to vanish if we take $A(w_1,w_2) \propto w_1-w_2$ and integrate both $w_1$ and $w_2$ around the point $z_1$. Let therefore
\be
 \psi _{12} = \Delta_t(z_1,\ldots,z_N)\oint \frac{1}{2\pi\i}\ \oint \frac{1}{2\pi\i}\ \phi_{12}(w_1,w_2,)\hspace{1mm}\d w_1\d w_2,
\ee
where the integration contour encircles the poles at $w=z_i$ but not those at $w=t^{-2}z_i$. In order to satisfy the initial condition
\be
 \psi _{12} = \Delta_t(z_1,z_2)\Delta_t(z_3,\ldots,z_N),
\ee
we find that
\be
A(w_1,w_2) = (w_2-w_1)(tw_1-t^{-1}w_2).
\ee

Finally we define

\be
\psi_{a_0a_1}(v_0,v_1) :=T_{a_0-1}(v_0+1) \cdots T_1(v_0+a_0-1)\cdot T_{a_1-2}(v_1+1) \cdots T_2(v_1+a_1-2)\ \psi_{12}
\ee

which we denote by 

\[
\psi_{a_0a_1}(v_0,v_1)= 
\begin{squarecells}{8}{3.5em}
$\scriptstyle v_1+a_1-2$ &  &  &  & & & $\scriptstyle v_1+2$ & $\scriptstyle v_1+1$ \nl
$\scriptstyle v_0+a_0-1$ & & & $\scriptstyle v_0+2$ & $\scriptstyle v_0+1$ \cr
\cline{1-5}
\end{squarecells}\ .
\]

\begin{lemma}

Let $y_k = -\frac{[v_k]}{[v_k+1]}$. For $a_1 > a_0 \geq 1$, \\ \\
\begin{align*}
\psi_{a_0a_1}(v_0,v_1) &=  \Delta_t(z_1,\ldots,z_N) \oint \frac{1}{2\pi\i}\oint \frac{1}{2\pi\i}\ (w_2-w_1) (tw_1-t^{-1}w_2) \times\\
& \hphantom{\Delta_t(z_1,\ldots,z_N)} \left(1+y_{ 0 } S_{a_{ 0 }}\right)\ \left(1+y_{ 1 } S_{a_{ 1 }}\right)\cdot \phi_{a_{ 0 }}(w_1)\phi_{a_{ 1 }}(w_2)\hspace{1mm}\d w_1\d w_2\\
&= \Delta_t(z_1,\ldots,z_N) \oint \frac{1}{2\pi\i}\oint \frac{1}{2\pi\i}\ (w_2-w_1) (tw_1-t^{-1}w_2) \times\\
&\hphantom{\Delta_t(z_1,\ldots,z_N)} \prod_{m=0}^1 \left(\frac{1}{[v_m+1]}\frac{t^{v_m+1} w_{m+1  } -t^{-v_m-1}z_{a_m}}{t w_{m+1  }-t^{-1}z_{a_m}} \right) \phi_{a_m}(w_{m+1  })\hspace{1mm}\d w_1\d w_2.
\end{align*}
\end{lemma}
\begin{proof}
The proof again uses the induction argument of Lemma~\ref{lem:onerow}, as well as the fact that $T_a(v)$ commutes with functions symmetric in $z_1,\ldots,\hat{z}_{a},\hat{z}_{a+1},\ldots, z_N$.
\end{proof}
\bigskip
\subsection{$m$ rows}
The proof of Theorem~\ref{DefMacdIntegral} can now be completed in similar fashion. Recall that we have $a_k=\lambda_k+k+1$ and $y_k = -\frac{[v_k]}{[v_k+1]}$. We have:
\begin{align*}
\psi_{a_0a_1\ldots a_{n-1}}(v_0,v_1,\ldots,v_{n-1}) &=\Delta_t(z_1,\ldots,z_N) \oint  \ldots \oint  \frac{1}{(2\pi\i)^n} \Delta(w_n,\ldots,w_1) \Delta_t(w_1,\ldots,w_n) \times\\
&\hphantom{\Delta_t(z_1)} \prod_{m=0}^{n-1} \left(1+y_m S_{a_m}\right)\ \phi_{a_m}(w_{m+1})\hspace{1mm}\d w_1 \dots \d w_{n}\\
&= \Delta_t(z_1,\ldots,z_N) \oint  \ldots \oint \frac{1}{(2\pi\i)^n}  \Delta(w_n,\ldots,w_1) \Delta_t(w_1,\ldots,w_{n}) \times\\
&\hphantom{\Delta_t(z_1)} \prod_{m=0}^{n-1} \frac{1}{[v_m+1]}\left(\frac{t^{v_m+1 } w_{m+1} -t^{-v_m-1}z_{a_m}}{t w_{m+1}-t^{-1}z_{a_m}} \right) \phi_{a_m}(w_{m+1})\hspace{1mm}\d w_1 \dots \d w_{n}.
\end{align*}

\begin{proof}
As before, by induction on $a_{0},\ldots,a_{n-1}$, and the fact that $T_a(v)$ commutes with functions symmetric in $z_1,\ldots,\hat{z}_{a},\hat{z}_{a+1},\ldots, z_N$.
\end{proof}
\bigskip
The polynomials $M_\lambda$ are now given by specialising $N=2n$, $a_0=1$ and setting $v_k=u_k+k-\lambda_k$, recalling that $v_0=0$ and hence $y_0=0$.    
\[
M_{(\lambda_{n-1},\ldots,\lambda_1)}(u_1,\ldots,u_{n-1}) = \psi_{1a_1\ldots a_{n-1}}(0,v_1,\ldots,v_{n-1})
\]
This is demonstrated in the following example: 
\bigskip
\begin{example}[$N=2n=6$]
Let 
\[
\M_{(1,1)}(u_1,u_2)=\psi_{1,3,4}(0,v_1,v_2)=
\raisebox{30pt}{
\begin{picture}(60,90)
\put(0,0){
\begin{squarecells}{1}{30pt}
$\scriptstyle u_2+2$ \cr\cline{1-1}
$\scriptstyle u_1+1$ \cr\cline{1-1}
\end{squarecells}}
\put(4,-28){\line(0,-1){30}}
\end{picture}}
= T_2(u_1+1) T_3(u_2+2) \Delta_t(z_1,z_2,z_3) \Delta_t(z_4,z_5,z_6).
\]
\vskip30pt

\noindent
 Since we always have $a_0=1$ we will from now refrain from drawing the bottom vertical line.  The integral representation for $\M_{(1,1)}(u_1,u_2)=\psi_{1,3,4}(0,v_1,v_2)$\ is
\begin{align*}
\M_{(1,1)}(u_1,u_2) &=  \Delta_t(z_1,\ldots,z_N) \oint \frac{1}{2\pi\i}\oint \frac{1}{2\pi\i}\oint \frac{1}{2\pi\i}\Delta(w_3,w_2,w_1)\Delta_t(w_1,w_2,w_3)\ \phi_1(w_1) \times\\
& \hphantom{\Delta_t(z_1,\ldots,z_N)} \left(1+y_1 S_{3}\right)\ \left(1+y_2 S_{4}\right)\cdot \phi_{3}(w_{2})\phi_{4}(w_{3})\hspace{1mm}\d w_1\d w_2\d w_3\\
&= \Delta_t(z_1,\ldots,z_N) \oint \frac{1}{2\pi\i}\oint \frac{1}{2\pi\i} \oint \frac{1}{2\pi\i} \Delta(w_3,w_2,w_1)\Delta_t(w_1,w_2,w_3)\ \phi_1(w_1) \times\\
& \hphantom{\Delta_t(z_1,\ldots,z_N)} \left(\frac{1}{[u_1+1]}\frac{t^{u_1+1} w_1 -t^{-u_1-1}z_{3}}{tw_1-t^{-1}z_{3}} \right) \phi_{3}(w_2)\times\\
& \hphantom{\Delta_t(z_1,\ldots,z_N)} \left(\frac{1}{[u_2+2]}\frac{t^{u_2+2} w_2 -t^{-u_2-2}z_{4}}{tw_2-t^{-1}z_{4}} \right) \phi_{4}(w_3)\hspace{1mm}\d w_1\d w_2 \d w_3.
\end{align*}
\end{example}
\section{Proof of Proposition~\ref{th:homlimit}}
\label{ap:homlimit}
Recalling that $\tau=-(t+t^{-1})$, the change of variables \eqref{eq:changeofvar} results in
\[
\widetilde{\Delta}_t(\zeta_1,\ldots,\zeta_{2n}) = \Delta_t(z_1,\ldots,z_{2n})= \frac{(t-t^{-1})^{n(2n-1)}}{\prod_{i=1}^{2n} (1-t\zeta_i)^{2n-1}} \prod_{1\leq i<j\leq 2n} (1+\zeta_i\zeta_j+\tau \zeta_j)
\]
and
\[
\phi_i(x) =  \prod_{j=1}^{2n} \frac{(1-t\zeta_j)(1-tx)}{t-t^{-1}} \prod_{j=1}^i \frac{1}{x-\zeta_j} \prod_{j=i+1}^{2n} \frac{1}{1+x\zeta_j+\tau \zeta_j}\ .
\]
In the homogeneous limit, $\zeta_j\rightarrow 0$, this reduces to
\[
\phi_i(x) \rightarrow \left(\frac{1-tx}{t-t^{-1}}\right)^{2n} \frac{1}{x^i}\ .
\]
We thus find that $\M_{(\lambda_{n-1},\ldots,\lambda_1)}(u_1,\ldots,u_{n-1})$\ can be written as 

\begin{multline*}
\M_{(\lambda_{n-1},\ldots,\lambda_1)}(u_1,\ldots,u_{n-1}) = (t-t^{-1})^{n(n-1)}\prod_{1\leq i<j\leq 2n}(1+\zeta_i\zeta_j+\tau \zeta_j) \prod_{j=1}^{2n} \frac{1}{(1-t \zeta_j)^{n-1}} \times\\
\oint  \cdots \oint \frac{1}{(2\pi\i)^n}\ \prod_{1\leq i<j\leq n} (x_j-x_i)(1+x_i x_j+\tau x_j) 
\prod_{m=0}^{n-1} \left(\frac{1+x_{m+1} \zeta_{a_m}- x_{m+1} \frac{[v_m]}{[v_m+1]}-\zeta_{a_m}\frac{[v_m+2]}{[v_m+1]}}{1+x_{m+1} \zeta_{a_m}+\tau \zeta_{a_m}} \right.\times\\
\left. \prod_{i=1}^{a_m} \frac{1}{x_{m+1}-\zeta_i} \prod_{i=a_m+1}^{2n} \frac{1}{1+x_{m+1}\zeta_i+\tau \zeta_i} \right)  \hspace{1mm}\d x_1 \dots \d x_n,
\end{multline*}

where $a_k=\lambda_k+k+1$, $v_k=u_k+k-\lambda_k$ and\ the integration is around the points $\zeta_i$. The homogenous limit $\zeta_i\rightarrow 0$ can now easily be taken.
\bigskip

\section{Expansion of $\M_{(3,2,1)}$ in terms of KL polynomials}
\label{ap:expansion}

We consider $n=4$.  We write the expansion \eqref{M1expansion} of the deformed maximal Macdonald polynomial $\M_{(2,1,0)}(u_2,u_3)$ as,
\ben
\begin{squarecells}{2}{2em}
$\scriptstyle 2+u_3$&$\scriptstyle 1+u_3$ \nl
$\scriptstyle 1+u_2$\cr
\cline{1-1}
\end{squarecells}
=
c_{(2,1)}\ \begin{squarecells}{2}{1em}
$\scriptstyle 2$&$\scriptstyle 1 $\nl
$\scriptstyle 1$\cr
\cline{1-1}
\end{squarecells}
+
c_{(1,1)}\ \begin{squarecells}{1}{1em}
$\scriptstyle 2$\nl
$\scriptstyle 1$\nl
\end{squarecells}
+
c_{(2)}\ \begin{squarecells}{2}{1em}
$\scriptstyle 2$&$\scriptstyle 1 $\nl
\end{squarecells}
+
c_{(1)}\ \begin{squarecells}{1}{1em}
$\scriptstyle 1$\nl
\end{squarecells}
+c_{\empty}\ \emptyset.
\een
Now, using \eqref{M2expansion} we can write the expansion for $\M_{(3,2)}$,
\begin{multline}
\label{eq:320exp}
M_{(3,2)}(u_2,u_3)=
\begin{squarecells}{3}{2em}
$\scriptstyle 3+u_3$&$\scriptstyle 2+u_3$&$\scriptstyle 1+u_3$ \nl
$\scriptstyle 2+u_2$&$\scriptstyle 1+u_2$ \cr \cline{1-2}
\end{squarecells}
=
c_{(2,1)}\ \begin{squarecells}{3}{1em}
$\scriptstyle 3$&$\scriptstyle 2$&$\scriptstyle 1 $\nl
$\scriptstyle 2$&$\scriptstyle 1$\cr
\cline{1-2}
\end{squarecells}
+
c_{(1,1)}\ \begin{squarecells}{2}{1em}
$\scriptstyle 3$&$\scriptstyle 2$\nl
$\scriptstyle 2$&$\scriptstyle 1$\nl
\end{squarecells}
+\\
c_{(2)}\ \begin{squarecells}{3}{1em}
$\scriptstyle 3$&$\scriptstyle 2$&$\scriptstyle 1 $\nl
$\scriptstyle 1$\cr \cline{1-1}
\end{squarecells}
+
c_{(1)}\ \begin{squarecells}{2}{1em}
$\scriptstyle 2$&$\scriptstyle 1$\nl
$\scriptstyle 1$ \cr \cline{1-1}
\end{squarecells}
+c_{\emptyset}\ \begin{squarecells}{1}{1em}
$\scriptstyle 2$\nl
$\scriptstyle 1$ \nl
\end{squarecells}.
\end{multline}
We want to obtain the expansion for $\M_{(3,2,1)}(u_1,u_2,u_3)$, and so by acting with $T_2(u_1+1)$ we obtain
\bea\label{KLexpandex}
&&M_{(3,2,1)}(u_1,u_2,u_3)=T_2(u_1+1) M_{(3,2)}(u_2,u_3)=(T_2(1)+y_1)M_{(3,2)}(u_2,u_3)\nn \\ &&\nn \\
&=&  
c_{(2,1)}\ \begin{squarecells}{3}{1em}
$\scriptstyle 3$&$\scriptstyle 2$&$\scriptstyle 1 $\nl
$\scriptstyle 2$&$\scriptstyle 1$\cr
\cline{1-2}
$\scriptstyle 1$\cr
\cline{1-1}
\end{squarecells}
+
c_{(1,1)}\ \begin{squarecells}{2}{1em}
$\scriptstyle 3$&$\scriptstyle 2$\nl
$\scriptstyle 2$&$\scriptstyle 1$\nl
$\scriptstyle 1$\cr
\cline{1-1}
\end{squarecells}
+
c_{(2)}\ \begin{squarecells}{3}{1em}
$\scriptstyle 3$&$\scriptstyle 2$&$\scriptstyle 1 $\nl
$\scriptstyle 1$\cr \cline{1-1}
$\scriptstyle 1$\cr
\cline{1-1}
\end{squarecells}
+
c_{(1)}\ \begin{squarecells}{2}{1em}
$\scriptstyle 2$&$\scriptstyle 1$\nl
$\scriptstyle 1$ \cr \cline{1-1}
$\scriptstyle 1$\cr
\cline{1-1}
\end{squarecells}
+c_{\emptyset}\ \begin{squarecells}{1}{1em}
$\scriptstyle 2$\nl
$\scriptstyle 1$ \nl
$\scriptstyle 1$\cr
\cline{1-1}
\end{squarecells} \nn
\\
&&+
y_1 \left( c_{(2,1)}\ \begin{squarecells}{3}{1em}
$\scriptstyle 3$&$\scriptstyle 2$&$\scriptstyle 1 $\nl
$\scriptstyle 2$&$\scriptstyle 1$\cr
\cline{1-2}
\end{squarecells}
+
c_{(1,1)}\ \begin{squarecells}{2}{1em}
$\scriptstyle 3$&$\scriptstyle 2$\nl
$\scriptstyle 2$&$\scriptstyle 1$\nl
\end{squarecells}
+
c_{(2)}\ \begin{squarecells}{3}{1em}
$\scriptstyle 3$&$\scriptstyle 2$&$\scriptstyle 1 $\nl
$\scriptstyle 1$\cr \cline{1-1}
\end{squarecells}
+
c_{(1)}\ \begin{squarecells}{2}{1em}
$\scriptstyle 2$&$\scriptstyle 1$\nl
$\scriptstyle 1$ \cr \cline{1-1}
\end{squarecells}
+c_{\emptyset}\ \begin{squarecells}{1}{1em}
$\scriptstyle 2$\nl
$\scriptstyle 1$ \nl
\end{squarecells} \hspace{2mm}\right).
\eea
Now, the third, fourth and fifth elements in this expansion are \underline{not} KL elements;  however, they may be expanded in terms of KL elements, for example we can write
\be\label{klexpand}
\begin{squarecells}{3}{1em}
$\scriptstyle 3$&$\scriptstyle 2$&$\scriptstyle 1 $\nl
$\scriptstyle 1$\cr \cline{1-1}
$\scriptstyle 1$\cr
\cline{1-1}
\end{squarecells}
= \hspace{2mm}
\begin{squarecells}{3}{1em}
$\scriptstyle 3$&$\scriptstyle 2$&$\scriptstyle 1 $\nl
$\scriptstyle 2$\cr \cline{1-1}
$\scriptstyle 1$\cr
\cline{1-1}
\end{squarecells}
+ \hspace{2mm}
\begin{squarecells}{3}{1em}
$\scriptstyle 3$&$\scriptstyle 2$&$\scriptstyle 1 $\nl
\end{squarecells}.
\ee 

Now we see how this result is obtained using the expansion over Dyck ribbons in \eqref{KLrecursion}, which we recall here:
\be
T_2(u_1+1)\ \KL_{\lambda}=\KL_{(\lambda,1)} + y_1 \KL_{\lambda} + \sum_{d_{\rm d} \in R_\rho} \KL_{\lambda - r_{\rm d}}. 
\label{KLrecursion2}
\ee
If we consider again the third element on the right hand side of \eqref{KLexpandex}, we note that it arises from the third element in the right hand side of \eqref{eq:320exp}, i.e. the partition $\lambda=(3,1)$.
We thus need to consider all Dyck ribbons with the diagonal as shown:
\begin{figure}[here!]
\begin{center}
\resizebox{1.3cm}{!}{\includegraphics{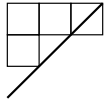}}
\caption{The diagonal corresponding to $\lambda=(3,1)$\label{diag310}}
\end{center}
\end{figure}

\begin{figure}[here!]
\begin{center}
\resizebox{1.3cm}{!}{\includegraphics{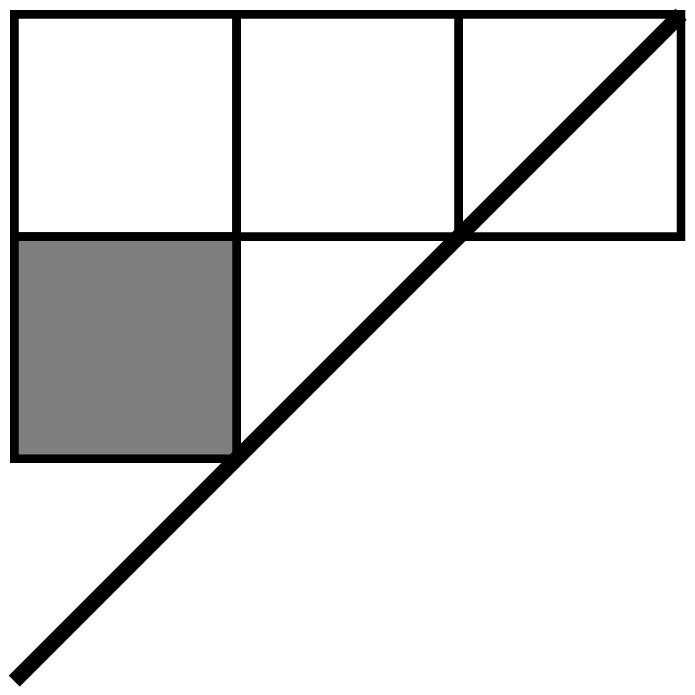}}
\caption{The Dyck ribbon
\label{ribb310}}
\end{center}
\end{figure}

There is only one such Dyck ribbon, shown shaded in Figure \ref{ribb310}. 
So the only partition $\lambda-r_{\rm d}$ appearing in the sum on the right hand side of \eqref{KLrecursion2} is $(3)$. Thus, from  \eqref{KLrecursion2}, the action of $T_2(u_1+1)$ on KL$_{(3,1)}$ is
\ben
T_2(u_1+1) \hspace{2mm}
\begin{squarecells}{3}{1em}
$\scriptstyle 3$&$\scriptstyle 2$&$\scriptstyle 1 $\nl
$\scriptstyle 1$\cr \cline{1-1}
\end{squarecells}
=
\begin{squarecells}{3}{1em}
$\scriptstyle 3$&$\scriptstyle 2$&$\scriptstyle 1 $\nl
$\scriptstyle 2$\cr \cline{1-1}
$\scriptstyle 1$\cr
\cline{1-1}
\end{squarecells}
+ \hspace{2mm}
y_1 \begin{squarecells}{3}{1em}
$\scriptstyle 3$&$\scriptstyle 2$&$\scriptstyle 1 $\nl
$\scriptstyle 1$\cr \cline{1-1}
\end{squarecells}
+ \hspace{2mm}
\begin{squarecells}{3}{1em}
$\scriptstyle 3$&$\scriptstyle 2$&$\scriptstyle 1 $\nl
\end{squarecells},
\een
in agreement with the result from the expansion  \eqref{KLexpandex}, using the \eqref{klexpand}.

From \eqref{KLexpandex} we can also deduce recursions between the coefficients $c_{\lambda}$. For example, it is clear that $c_{3,2,1}=c_{2,1}$, and that $c_{3,2}=y_1 c{2,1}$. In fact, all coefficients corresponding to two-row diagrams have a factor $y_1$, while others do not contain such a factor. Careful considerations lead to the recursion as explained in the main text.

\newpage

\newcommand\arxiv[1]{\href{http://arxiv.org/abs/#1}{\tt arXiv:#1}}
\bibliographystyle{h-model1-num-names}
\bibliography{poly}

\end{document}